\documentclass[11pt]{article}

\usepackage{amssymb,amsmath,amsthm}
\usepackage{color}
\usepackage{pdflscape}
\usepackage{longtable}
\usepackage{enumerate}
\usepackage{tikz}
\usepackage{float}
\usepackage[colorlinks=true, linkcolor=blue, urlcolor=blue, citecolor=blue]{hyperref}
\usepackage{url,hyperref}

\usepackage[margin=2.5cm]{geometry}

\usepackage[toc,page]{appendix}
\usepackage{lscape}
\usepackage{multirow}
\usepackage{adjustbox,expl3,etoolbox}
\usepackage{color} 
\usepackage[utf8]{inputenc}

\newtheorem{thm}{Theorem}[section]
\newtheorem{lem}[thm]{Lemma}
\newtheorem{prop}[thm]{Proposition}
\theoremstyle{definition}
\newtheorem{defn}[thm]{Definition}
\newtheorem{conj}[thm]{Conjecture}
\newtheorem{cor}[thm]{Corollary}

\newtheorem{rmk}[thm]{Remark}
\newtheorem{prob}[thm]{Problem}
\newtheorem{qst}[thm]{Question}

\newtheorem{hyp}{Hypothesis}

\newcommand{\itbf}[1]{{\bf{{\emph{{#1}}}}}}
\newcommand{\Aut}[1]{\operatorname{Aut}(#1)}
\newcommand{\End}[2]{\operatorname{End}_{#1}\hspace*{-0.1cm}\left(#2\right)}
\newcommand{\sym}[1]{\operatorname{Sym}(#1)}
\newcommand{\alt}[1]{\operatorname{Alt}(#1)}
\newcommand{\fis}{\operatorname{Fi}}

\newcommand{\agl}[2]{\operatorname{AGL}_#1(#2)}

\newcommand{\agammal}[2]{\operatorname{A\Gamma L}_#1(#2)}

\newcommand{\mathieu}[1]{\operatorname{M}_{#1}}

\newcommand{\sln}[2]{\operatorname{SL}_{#1}(#2)}
\newcommand{\pg}[2]{\operatorname{PG}_{#1}(#2)}

\newcommand{\qbinom}[3]{\begin{bmatrix}
		#1\\ #2
	\end{bmatrix}_{#3}}

\usepackage{setspace,kantlipsum}
\newcommand{\vo}{\mathrm{VO}}
\newcommand{\vd}{\mathrm{VD}}
\newcommand{\vsz}{\mathrm{VSz}}
\newcommand{\no}{\mathrm{NO}}
\newcommand{\nug}{\mathrm{NU}}
\newcommand{\vLS}{\operatorname{vLS}}

\makeatletter
\renewcommand*{\@fnsymbol}[1]{\@alph{#1}}
\makeatother

\title{\bf On the classification of triply-transitive strongly-regular graphs}

\author{
Allen Herman\footnote{The first author's research is supported by a grant from NSERC.}\\
\small Department of Mathematics and Statistics \\[-0.8ex]
\small University of Regina \\[-0.8ex] 
\small 3737 Wacana Parkway, Regina, Saskatchewan, Canada \\
\small\tt {allen.herman}@uregina.ca  \vspace*{0.5cm}\\	
Roghayeh Maleki\thanks{Corresponding author: \href{mailto:roghayeh.maleki@uregina.ca}{roghayeh.maleki@uregina.ca}} 
\\
\small Department of Mathematics and Statistics \\[-0.8ex]
\small University of Regina \\[-0.8ex] 
\small 3737 Wacana Parkway, Regina, Saskatchewan, Canada \\
\small\tt {roghayeh.maleki}@uregina.ca 
\vspace*{0.5cm}\\
 Andriaherimanana Sarobidy Razafimahatratra\footnote{The third author gratefully acknowledges that this research was supported by the Fields Institute for Research
 	in Mathematical Sciences.}\\
\small Fields Institute for Research in Mathematical Sciences, \\[-0.8ex]
\small 222 College St., Toronto, ON M5T 3J1, Canada \\
\small\tt sarobidy@phystech.edu \\
}

\begin{document}	

\maketitle

\begin{abstract}
 Let $\Gamma = (\Omega,E)$ be a strongly-regular graph with adjacency matrix $A_1$, and let $A_2$ be the adjacency matrix of its complement. For any vertex $\omega\in \Omega$, we define $E_{0,\omega}^*$ $E_{1,\omega}^*$ and $E_{2,\omega}^*$ to be respectively the diagonal matrices whose main diagonal is the row corresponding to $\omega$ in the matrices $I, A_1$, and $A_2$. The Terwilliger algebra of $\Gamma$ with respect to the vertex $\omega\in \Omega$ is the subalgebra $T_\omega = \left\langle I,A_1,A_2,E_{0,\omega}^*,E_{1,\omega}^*,E_{2,\omega}^* \right\rangle$ of the complex matrix algebra $\operatorname{M_{|\Omega|}}(\mathbb{C})$. The algebra $T_\omega$ contains the subspace $T_{0,\omega} = \operatorname{Span}\left\{ E_{i,\omega}^*A_jE_{k,\omega}^*: 0\leq i,j,k\leq 2 \right\}$.
 In addition, if $G = \Aut{\Gamma}$, then $T_\omega$ is a subalgebra of the centralizer algebra $\tilde{T}_\omega = \End{G_\omega}{\mathbb{C}^\Omega}$.
 The strongly-regular graph $\Gamma=(\Omega,E)$ is triply transitive if $\Gamma$ is vertex transitive and $T_{0,\omega} = T_\omega = \tilde{T}_\omega$, for any $\omega \in \Omega$. In this paper, we classify all triply transitive strongly-regular graphs that are not isomorphic to the collinearity graph of the polar space $O_{6}^-(q)$, where $q$ is a prime power, or the affine polar graph $\vo_{2m}^\varepsilon(2)$, where $m\geq 1$ and  $\varepsilon = \pm 1$.
  \\ 
 \noindent {\bf Mathematics Subject Classifications:} 05E30; 20B05, 51E99.
 
 \noindent {\bf Keywords:} Strongly-regular graphs, Terwilliger algebras, automorphism groups, finite geometry.
\end{abstract}

\section{Introduction}
Let $\Omega$ be a finite set and $\mathcal{R} = \{ R_0,R_1,\ldots,R_d \}$ be a partition of $\Omega \times \Omega$. The pair $(\Omega,\mathcal{R})$ is called an \itbf{association scheme} if the following properties are satisfied
\begin{enumerate}
	\item $R_0 = \{ (\omega,\omega): \omega\in \Omega \}$,
	\item  for all $1\leq i\leq d$, the relation $R_{i^*} = \{ (\omega^\prime,\omega): (\omega,\omega^\prime) \in R_i \}$ belongs to $\mathcal{R}$,
	\item for any $0\leq i,j,k\leq d$, there exists a non-negative integer $\mathrm{p}_{ij}^k$ such that for $(\omega,\omega^\prime) \in R_k$, we have
	\begin{align*}
		\mathrm{p}_{ij}^k = \left|\left\{ \delta \in \Omega: (\omega,\delta) \in R_i, (\delta,\omega^\prime) \in R_j \right\}\right|,
	\end{align*}
	\item $\mathrm{p}_{ij}^k = \mathrm{p}_{ji}^k$, for any $0\leq i,j,k\leq d$.
\end{enumerate}
In this case, $(\Omega,\mathcal{R})$ is called a $d$-class association scheme.

Let $(\Omega,\mathcal{R})$ be an association scheme, where $\mathcal{R} = \{R_0,R_1,\ldots,R_d\}$. The pair $(\Omega,R_i)$, for any $1\leq i\leq d$, is a digraph on $\Omega$ and arc set $R_i$. For any $1\leq i\leq d$, let $A_i$ be the adjacency matrix of the digraph $(\Omega,R_i)$, and assume that $A_0 = I_{|\Omega|}$ is the identity matrix of order $|\Omega|$. The {Bose-Mesner algebra $\mathcal{B}_\mathcal{R}$} of the association scheme $(\Omega,\mathcal{R})$ is the complex subalgebra of the matrix algebra $\operatorname{M_{|\Omega|}}(\mathbb{C})$ generated by the adjacency matrices $A_0,A_1,\ldots,A_d$. Since 
\begin{align*}
	A_iA_j = \sum_{k=0}^{d}\mathrm{p}_{ij}^k A_k
\end{align*}
for any $0\leq i,j\leq d$, the Bose-Mesner algebra is also spanned by $\{A_0,A_1,\ldots,A_d\}$. The Bose-Mesner algebra of the association scheme $(\Omega,\mathcal{R})$ acts naturally by left multiplication on the vectors of $\mathbb{C}^{|\Omega|}$. This action endows $\mathbb{C}^{|\Omega|}$ with a $\mathcal{B}_\mathcal{R}$-module structure. We will call the latter the \itbf{standard module} of the association scheme $(\Omega,\mathcal{R}).$

Now fix $\omega\in \Omega$. For any $0\leq i\leq d$, define the matrix $E_{i,\omega}^*$, whose rows and columns are indexed by $\Omega$, to be the diagonal matrix with entries
\begin{align}
	E^*_{i,\omega}(\alpha,\alpha) =
	\begin{cases}
		1 &\mbox{ if } (\omega,\alpha) \in R_i,\\
		0 &\mbox{ otherwise}.
	\end{cases} 
	\label{eq:idempotents}
\end{align}
The \itbf{Terwilliger algebra} of the association scheme $(\Omega,\mathcal{R})$ with respect to $\omega \in \Omega$ is the algebra
\begin{align}
	T_{\omega}((\Omega,\mathcal{R})) = \left\langle A_0,A_1,\ldots,A_d, E_{0,\omega}^*,E_{1,\omega}^*,\ldots,E_{d,\omega}^* \right\rangle.
\end{align}
These algebras were introduced by Paul Terwilliger as subconstituent algebras in \cite{terwilliger1992subconstituent,terwilliger1993subconstituent2,terwilliger1993subconstituent3}. We will usually just write $T_\omega$ instead of $T_{\omega}((\Omega,\mathcal{R}))$, if the association scheme $(\Omega,\mathcal{R})$ is clear from the context. The standard module of the underlying association scheme also admits a Terwilliger algebra module structure.

\subsection{Motivation}
The automorphism group $\Aut{\mathcal{R}}$ of the association scheme $(\Omega,\mathcal{R})$ is the set of all permutations of $\Omega$ that leave every relation in $\mathcal{R}$ invariant. For any permutation group $H\leq \sym{\Omega}$, the centralizer algebra $\End{H}{\mathbb{C}^{|\Omega|}}$ is the set of all matrices in $\operatorname{M_{|\Omega|}}(\mathbb{C})$ that commute with all elements of $H$ viewed as permutation matrices on $\Omega$. For any $\omega \in \Omega$, we always have the inclusion
\begin{align}
	T_\omega \subseteq \End{\Aut{\mathcal{R}}_\omega}{\mathbb{C}^{|\Omega|}},\label{eq:inclusion1}
\end{align}
see \cite[Lemma~2.2]{HMR2025}.
If $\Aut{\mathcal{R}}$ acts transitively on $\Omega$, then it is not hard to show that for any $\omega,\omega^\prime \in \Omega$, the Terwilliger algebra $T_\omega$ and $T_{\omega^\prime}$ are isomorphic. Hence, if $\Aut{\mathcal{R}}$ is transitive on $\Omega$, then we will just denote the Terwilliger algebra with respect to any element of $\Omega$ by $T$. Similarly, we denoted the centralizer algebra of the point stabilizer $\Aut{\mathcal{R}}_\omega$ by $\tilde{T}$ instead of $\End{\Aut{\mathcal{R}}_\omega}{\mathbb{C}^{|\Omega|}}$. In this case, \eqref{eq:inclusion1} becomes $T\subseteq \tilde{T}$.

We note that the gap in dimension between the algebras in \eqref{eq:inclusion1} can be arbitrarily large, as in the case of Paley graphs on a prime number of vertices \cite[Corollary~8.6]{hanaki2023terwilliger}. However, equality holds in \eqref{eq:inclusion1} for many Terwilliger algebras of classical association schemes arising from distance-regular graphs such as the Johnson graphs \cite{tan2019terwilliger}, Hamming graphs \cite{gijswijt2006new,schrijver2005new}, and folded $n$-cubes \cite{hou2020new}. Moreover, equality also holds in \eqref{eq:inclusion1} for association schemes such as conjugacy class schemes of abelian groups \cite{bannai1995terwilliger} and metacyclic groups of the form $C_n\rtimes C_2$ \cite{maleki2024terwilliger}. 

The existence of these examples where equality holds in \eqref{eq:inclusion1} is not too surprising. In fact, in the very first paper \cite{terwilliger1992subconstituent} where Terwilliger algebras were introduced, Paul Terwilliger noted that $T_\omega$ is a combinatorial analogue of the centralizer algebra $\End{\Aut{\mathcal{R}}_\omega}{\mathbb{C}^{|\Omega|}}$. The centralizer algebra of $\Aut{\mathcal{R}}_\omega$ has a rich algebraic structure which is completely understood from the representation theory of the underlying group, whereas the Terwilliger algebra is sometimes much more complicated due to the underlying combinatorial structures. Therefore, association schemes for which equality holds in \eqref{eq:inclusion1} simultaneously enjoy rich algebraic and combinatorial properties. Due to this reason, they are of interests to many researchers.

In \cite{tan2019terwilliger}, the authors posed a question essentially asking for which association schemes do we have equality in \eqref{eq:inclusion1}. We state this in the next problem.
\begin{prob}
	Given an association scheme $(\Omega,\mathcal{R})$, when does equality hold in \eqref{eq:inclusion1}?\label{prob}
\end{prob}

The above-mentioned problem is too general so one cannot expect a comprehensive answer. Instead, we will focus on studying Problem~\ref{prob} for association schemes with few classes. Strongly-regular graphs are the first natural examples to consider since they correspond to symmetric $2$-class association schemes, implying that the automorphism group of the scheme is equal to the automorphism group of any of the graphs in the scheme. This is crucial for \eqref{eq:inclusion1} since the automorphism group of an association scheme can be much smaller than the automorphism group of any graph in the scheme. Therefore, we will focus on this problem for association schemes corresponding to strongly-regular graphs. Even for strongly-regular graphs, Problem~\ref{prob} still seems to be intractable. So, we study the strongly-regular graphs where \eqref{eq:inclusion1} holds and admits an additional property called triple regularity.

\subsection{Main results}
Recall that a graph $\Gamma = (\Omega,E)$ is a \itbf{strongly-regular graph} with parameters $(v,k,\lambda,\mu)$ if $\Gamma$ is a $k$-regular graph on $v = |\Omega|$ vertices, and with the properties that any $\omega,\omega^\prime \in \Omega$ such that $\{\omega,\omega^\prime\} \in E$ are adjacent to exactly $\lambda$ vertices, and any $\omega,\omega^\prime \in \Omega$ such that $\{\omega,\omega^\prime\} \not \in E$ are adjacent to exactly $\mu$ vertices.
The rank of a transitive permutation group is the number of orbits of its point stabilizers, or equivalently, its number of orbitals. The automorphism group of a vertex-transitive strongly-regular graph has rank at least $3$. For example, the Shrikhande graph is a strongly-regular graph with full automorphism group of rank $4$. If the rank of the automorphism group is equal to $3$, then the strongly-regular graph is called a \itbf{rank $3$ graph.}

We shall assume the following hypothesis for most of this paper, so we state it here for convenience.

\begin{hyp}
	Let $\Gamma = (\Omega, E)$ be a strongly-regular graph with parameters $(n,k,\lambda,\mu)$. Let $\omega \in \Omega$ and define $\Delta_0 = \{\omega\}$, $\Delta_1$ the neighbours of $\omega$ in $\Gamma$, and $\Delta_2$ the neighbours of $\omega$ in $\overline{\Gamma}$. Let $\left(E_{i,\omega}^*\right)_{i\in \{0,1,2\}}$ be the idempotents defined in \eqref{eq:idempotents}, and consider the Terwilliger algebra $T_\omega$ with respect to $\omega \in \Omega$.\label{hyp0}
\end{hyp} 

Let $\Gamma = (\Omega,E)$ be a strongly-regular graph as in Hypothesis~\ref{hyp0}, and denote the automorphism group of $\Gamma$ by $G = \Aut{\Gamma}$. 
Define the subspace $T_{0,\omega}$ of $T_\omega$ to be 
\begin{align}
	T_{0,\omega} = \operatorname{Span}\left\{ E_{i,\omega}^*A_jE_{k,\omega}^*: 0\leq i,j,k\leq d \right\}.\label{eq:T0}
\end{align}
The graph $\Gamma$ is \itbf{triply regular} if for any $\alpha,\beta,\gamma\in \Omega$, the integer
\begin{align*}
	|\left\{ \omega\in \Omega: d(\alpha,\omega) = i,d(\beta,\omega) = j,d(\gamma,\omega)=k \right\}|
\end{align*}
is a constant depending only on $i,j,k$ and $a = d(\alpha,\beta), b =  d(\beta,\gamma),c = d(\alpha,\gamma)$.
In an unpublished note \cite[Lemma~4]{Munemasa}, Munemasa showed that the strongly-regular graph $\Gamma$ is triply-regular if and only if there exists an $\omega\in \Omega$ such that $T_{0,\omega} = T_\omega$. 

If the automorphism group of the strongly-regular graph $\Gamma$ is transitive, then we write $T_0$ instead of $T_{0,\omega}$. Next, we define the notion of triple transitivity.
\begin{defn}
	We say that the graph $\Gamma = (\Omega,E)$ is \itbf{triply transitive} if
	\begin{enumerate}[(1)]
		\item $\Aut{\Gamma}$ is transitive on $\Omega$, and
		\item $T_0 = T = \tilde{T}$.
	\end{enumerate} \label{def:1}
\end{defn}
Hence, if $\Gamma$ is triply transitive, then it is triply regular and \eqref{eq:inclusion1} holds with equality. In \cite{bannai1995terwilliger}, Bannai and Munemasa used a similar terminology when the association scheme is a conjugacy class scheme of an abstract group. However, the centralizer algebra $\tilde{T}$ in \cite{bannai1995terwilliger} is larger than $\End{\Aut{\mathcal{R}}_\omega}{\mathbb{C}^{|\Omega|}}$ since they did not consider a full point stabilizer.

Before stating the main result, we will need a few additional notions. The \itbf{Latin square graph} $\operatorname{LS}_m(n)$ is a strongly-regular graph arising from a transversal design $TD(m,n)$, where $2\leq m\leq n$ are integers. See \cite[\S8.4.2]{brouwer2022strongly} for the construction. The graph $\operatorname{LS}_m(n)$ has parameters
\begin{align*}
	\left(n^2,m(n-1),(m-1)(m-2)+n-2,m(m-1)\right).
\end{align*}
A strongly-regular graph with the same parameters as a Latin square graph but does not arise from a transversal design is said to have \itbf{Latin square parameters}. Such graphs are also known as \itbf{pseudo Latin square graphs}.
A strongly-regular graph $\Gamma$ has \itbf{negative Latin square parameters} if there exist non-negative integers $m$ and $n$ with $2\leq m\leq n$ such that the parameters of $\Gamma$ are
\begin{align*}
	\left(n^2,m(n+1),m(m+3)-n,m(m+1)\right).
\end{align*}

We state the first main result of the paper as follows.
\begin{thm}
	Let $\Gamma = (\Omega,E)$ be a strongly-regular graph as in Hypothesis~\ref{hyp0}. In addition, assume that the parameters of $\Gamma$ are  not Latin square or a negative Latin square parameters. Then, $\Gamma$ is triply transitive, then it is one of the following.
	\begin{enumerate}[(i)]
		\item A complete multipartite graph with $n$ parts of size $m$, where $n\neq m$.\label{imprimitive1}
		\item The cycle on $5$ vertices\footnote{this is also the Paley graph $P(5)$}.
		\item The McLaughlin graph.\label{McL1}
		\item The collinearity graph of the polar space  $O_6^-(q)$.\label{o}
	\end{enumerate}
	Conversely, all graphs in \eqref{imprimitive1}-\eqref{McL1} are triply transitive.
	\label{thm:main}
\end{thm}

The case \eqref{o} in Theorem~\ref{thm:main} is still open. However, they are triply transitive for small values of $q$.

The next results narrow down the families from which a triply transitive graph with Latin square or negative Latin square parameters can arise from.
\begin{thm}
	Let $\Gamma = (\Omega,E)$ be a strongly-regular graph as in Hypothesis~\ref{hyp0} and assume that $\Gamma$ has Latin square or a negative Latin square parameters. If $\Gamma$ is triply transitive, then it is one of the following families.
	\begin{enumerate}[(i)]
		\item The Higman-Sims graph.\label{HS}
		\item A complete multipartite graph with $n$ parts of size $n$, where $n\geq 2$ is an integer.\label{cmt}
		\item The Peisert graph $P^*(9)$ (which is isomorphic to the Paley graph on $9$ vertices).
		\item The $n\times n$ grid for some integer $n\geq 2$.\label{grid}
		\item The affine polar graph $VO^\varepsilon_{2m}(2)$, for $m\geq 2$ and $\varepsilon =\pm 1$.\label{vo}
	\end{enumerate}
	Conversely, the graphs in \eqref{HS}-\eqref{grid} are triply transitive.
	
	\label{thm:main2}
\end{thm}
The case \eqref{vo} in Theorem~\ref{thm:main2} is still open. They are however triply transitive for small values of $q$.

\subsection{Organization of the paper}
In Section~\ref{Sec:background}, we recall the necessary background results that are needed to prove the main results. In Section~\ref{Sec:TildeT}, we show that the point stabilizer of the full automorphism group of a triply transitive graph is a rank $3$ group. From this, we deduce some important results on its centralizer algebra. We then recall the classification of rank $3$ graphs in Section~\ref{Sec:rank3-graphs}. In Section~\ref{Sec:first-cases} and Section~\ref{subsect:general}, a comprehensive analysis of the triply transitive graphs is given, using the classification of rank $3$ graphs. The proofs of the main results, Theorem~\ref{thm:main} and Theorem~\ref{thm:main2}, follow immediately from this analysis, and they are given in Section~\ref{Sec:main-results}. In Section~\ref{Sec:conclusion}, we pose some open problems and make some conjectures about triply transitive graphs.

\section{Background results}\label{Sec:background}
Let $\Gamma = (\Omega,E)$ be a graph on $n$ vertices. Recall that $\Gamma$ is strongly regular with parameters $(n,k,\lambda,\mu)$ if 
\begin{enumerate}
	\item $\Gamma$ is $k$-regular,
	\item any $\alpha,\beta\in \Omega$ that are adjacent in $\Gamma$ have exactly $\lambda$ common neighbours,
	\item any $\alpha,\beta\in \Omega$ that are non-adjacent in $\Gamma$ have exactly $\mu$ neighbours.
\end{enumerate}
Strongly-regular graphs have been studied extensively in algebraic graph theory (see \cite{brouwer2022strongly}). If $\Gamma$ is strongly regular with parameters $(n,k,\lambda,\mu)$, then its complement $\overline{\Gamma}$ is also strongly regular with parameters $(n,n-k-1,n-2k+\mu-2,n-2k+\lambda)$.

Let $\Gamma = (\Omega,E)$ be a strongly-regular graph. Define
\begin{align}
	\begin{split}
		R_0 &= \{ (\omega,\omega):\omega \in \Omega \},\\ R_1 &= \{ (\omega,\omega^\prime): \{ \omega,\omega^\prime \}\in E \},\mbox{ and }\\ R_2 &= \left\{ (\omega,\omega^\prime): \{\omega,\omega^\prime\}\not\in E  \right\}.
	\end{split}\label{eq:relations}
\end{align}
The pair $\left(\Omega,\{R_0,R_1,R_2\}\right)$ is a $2$-class association scheme, and its intersection numbers are given by 
\begin{align}
	\begin{array}{lll}
		\mathrm{p}_{00}^0=1 & \mathrm{p}_{01}^0=0 & \mathrm{p}_{02}^0=0\\
		\mathrm{p}_{10}^1=1 & \mathrm{p}_{11}^1=\lambda & \mathrm{p}_{12}^1=k-\lambda-1\\
		\mathrm{p}_{20}^2=1 & \mathrm{p}_{21}^2=k-\mu & \mathrm{p}_{22}^2=n-2k+\mu-2\\
		\mathrm{p}_{00}^1=0 & \mathrm{p}_{01}^1=1 & \mathrm{p}_{02}^1=0\\
		\mathrm{p}_{00}^2=0 & \mathrm{p}_{01}^2=0 & \mathrm{p}_{02}^2=1\\
		\mathrm{p}_{10}^0=0 & \mathrm{p}_{11}^0=k & \mathrm{p}_{12}^0=0\\
		\mathrm{p}_{10}^2=0 & \mathrm{p}_{11}^2=\mu & \mathrm{p}_{12}^2=k-\mu\\
		\mathrm{p}_{20}^0=0 & \mathrm{p}_{21}^0=0 & \mathrm{p}_{22}^0=n-k-1\\
		\mathrm{p}_{20}^1=0 & \mathrm{p}_{21}^1=k-\lambda-1 & \mathrm{p}_{22}^1=n-2k+\lambda.\\
	\end{array}
	\label{eq}
\end{align}

The graph $\Gamma$ is called \itbf{primitive} if $\Gamma$ and its complement $\overline{\Gamma}$ are connected, otherwise it is called \itbf{imprimitive}. From \cite[Section~1.1.3]{brouwer2022strongly}, $\Gamma$ is imprimitive if and only if the intersection numbers given in \eqref{eq} satisfy $\mathrm{p}_{21}^1 = \mathrm{p}_{12}^1 = 0$ or $\mathrm{p}_{21}^2 = \mathrm{p}_{12}^2 = 0$, or equivalently, $\mathrm{p}_{11}^2 = 0$ or $\mathrm{p}_{22}^1 = 0$. If $\Gamma$ is triangle free, then $\mathrm{p}_{11}^1 = \lambda = 0$, and similarly, $\mathrm{p}_{22}^2 = n-2k+\mu - 2 = 0$ if $\overline{\Gamma}$ is triangle free. 
It is an easy exercise to verify that the intersection numbers in \eqref{eq} which depend on the parameters $(n,k,\lambda,\mu)$ are non-zero, unless one of the graphs in the association scheme is imprimitive or triangle free. 

\subsection{Terwilliger algebras of strongly-regular graphs}
Let $\Gamma = (\Omega,E)$ be a strongly-regular graph as in Hypothesis~\ref{hyp0}.
In \cite{tomiyama1994subconstituent}, Tomiyama and Yamazaki gave a complete description of the dimension of the Terwilliger algebra of strongly-regular graphs.
Their results rely on the subconstituents of $\Gamma$. By Hypothesis~\ref{hyp0}, we can partition the vertex set of $\Gamma$ into $\Delta_0=\{\omega\}$, $\Delta_1 = \{ \alpha\in \Omega: \{\omega,\alpha\} \in E \}$, and  $\Delta_2 = \{ \alpha\in \Omega: \{\omega,\alpha\} \not\in E \}$. The first and second subconstituents of $\Gamma$ are respectively the subgraphs $\Gamma_1$ and $\Gamma_2$ of $\Gamma$ induced by $\Delta_1$ and $\Delta_2.$ Let $B_1$ and $B_2$ be the adjacency matrices of $\Gamma_1$ and $\Gamma_2$, respectively. It is well known that the first and second subconstituents are regular graphs. Let $\tau<\theta<k$ be the eigenvalues of $\Gamma$. 

\begin{thm}
	\cite[Theorem~1.1]{tomiyama1994subconstituent}
	Let $\Gamma$ be a strongly-regular graph as in Hypothesis~\ref{hyp0}. The dimension of the Terwilliger algebra of $T_\omega$ of $\Gamma$ is equal to
	\begin{align*}
		m_1+m_2+4n_1+9 = m_1+m_2+4n_2+9,
	\end{align*}
	where 
	\begin{itemize}
		\item $m_1$ is the number of distinct eigenvalues of $\Gamma_1$ contained in $\{\tau,\theta\}$, corresponding to some eigenvectors orthogonal to the all-ones vector of size $k$,
		\item $m_2$ is the number of distinct eigenvalues of $\Gamma_2$ contained in $\{\tau,\theta\}$, corresponding to some eigenvectors orthogonal to the all-ones vector of size $n-k-1$,
		\item $n_1$ is the number of distinct eigenvalues of $\Gamma_1$ not contained in $\{\tau,\theta\}$, corresponding to some eigenvectors orthogonal to the all-ones vector of size $k$,
		\item $n_2$ is the number of distinct eigenvalues of $\Gamma_2$ not contained in $\{\tau,\theta\}$, corresponding to some eigenvectors orthogonal to the all-ones vector of size $n-k-1$.
	\end{itemize}\label{thm:dim-of-T}
\end{thm}

The result in Theorem~\ref{thm:dim-of-T} depends heavily on the knowledge of the spectra of the subconstituents of the underlying strongly regular graph. Since these subconstituents are just regular graphs, determining their eigenvalues might be highly non-trivial. Consequently, Theorem~\ref{thm:dim-of-T} is not practical, and we employ different techniques to study the Terwilliger algebra.
\subsection{Block decomposition of subspaces}
Let $\Gamma = (\Omega,E)$ be a strongly-regular graph as in Hypothesis~\ref{hyp0}.
Let $G = \Aut{\Gamma}$ and $\tilde{T}_\omega$ be the centralizer algebra of the point stabilizer $G_\omega = \{ g\in G: \omega^g = \omega \}$. 
We state the following proposition whose proof is straightforward.
\begin{prop}
	Let $U$ be a subspace of $\tilde{T}_\omega$ containing $(E_{i,\omega}^*)_{i=0,1,2}$. For any $0\leq i,j\leq 2$, the set $E_{i,\omega}^* U E_{j,\omega}^*$ is a subspace of $\tilde{T}_\omega$ and  $E_{i,\omega}^* \tilde{T}_\omega E_{j,\omega}^*$. If $U$ is a subalgebra of $\tilde{T}_\omega$, then so is $E_{i,\omega}^* U E_{i,\omega}^*$, for  $0\leq i\leq 2$.\label{prop:block-decomposition}
\end{prop}

The following corollary is also immediate.
\begin{cor}
	Let $U$ be a subspace of $\tilde{T}_\omega$ containing $(E_{i,\omega}^*)_{i=0,1,2}$. If $E_{i,\omega}^* U E_{i,\omega}^*$ is not closed under multiplication for some $0\leq i\leq 2$, then $U$ is not a subalgebra of $\tilde{T}_\omega$.
\end{cor}

Given a subspace $U$ of $\tilde{T}_\omega$ containing the mutually orthogonal idempotents $(E_{i,\omega}^*)_{i=0,1,2}$, the \itbf{Peirce decomposition} of $U$ is the direct sum of the subspaces $E^*_{i,\omega} U E_{j,\omega}^*$, for all $0\leq i,j\leq 2$.
The \itbf{block dimension decomposition} of $U$ is the matrix
\begin{align}
	\begin{bmatrix}
		\dim(E_{0,\omega}^* U E_{0,\omega}^*) & \dim(E_{0,\omega}^* U E_{1,\omega}^*) & \dim(E_{0,\omega}^* U E_{2,\omega}^*)\\
		\dim(E_{1,\omega}^* U E_{0,\omega}^*) & \dim(E_{1,\omega}^* U E_{1,\omega}^*) & \dim(E_{1,\omega}^* U E_{2,\omega}^*)\\
		\dim(E_{2,\omega}^* U E_{0,\omega}^*) & \dim(E_{2,\omega}^* U E_{1,\omega}^*) & \dim(E_{2,\omega}^* U E_{2,\omega}^*)\\
	\end{bmatrix}.
\end{align}
In particular, we have
\begin{align*}
	\dim(U) = \sum_{0\leq i,j\leq 2} \dim(E_{i,\omega}^* U E_{j,\omega}^*).
\end{align*}

Next, we determine the block dimension decomposition of $T_0$.	We recall the following result.
\begin{lem}[\cite{bannai1995terwilliger}]
	For any $0\leq i,k\leq 2$, we have
	\begin{align*}
		\dim(E_{i,\omega}^* T_0 E_{k,\omega}^*) = \left|\left\{  0\leq j\leq 2: \mathrm{p}_{ij}^k\neq 0\right\}\right|.
	\end{align*}
	 In particular, the dimension of the subspace $T_0 $ is equal to the number of non-zero intersection numbers $(\mathrm{p}_{ij}^k)$.
	\label{lem:dim-of-t0}
\end{lem}
By counting the non-zero intersection numbers in \eqref{eq}, the following proposition follows immediately.

\begin{prop}
	The possible block dimension decompositions of $T_0$ are given as follows.
	\begin{enumerate}[(i)]
		\item If $\Gamma$ is imprimitive, then the block dimension decomposition of $T_0$ is 
		\begin{align*}
			\begin{bmatrix}
				1 & 1 & 1\\
				1 & 2 & 1\\
				1 & 1 & 3
			\end{bmatrix}
			\mbox{ or }
			\begin{bmatrix}
				1 & 1 & 1\\
				1 & 2 & 1\\
				1 & 1 & 2
			\end{bmatrix}.
		\end{align*}
		In particular, $\dim(T_0) = 11$ and $\Gamma$ is a complete bipartite graph, or $\dim(T_0) = 12$ and $\Gamma$ is a complete multipartite graph with at least $3$ parts.\label{first}
		\item If $\Gamma$ is primitive, and both of $\Gamma$ and $\overline{\Gamma}$ are triangle-free, then  the block dimension decomposition of $T_0$ is 
		\begin{align*}
			\begin{bmatrix}
				1 & 1 & 1\\
				1 & 2 & 2\\
				1 & 2 & 2
			\end{bmatrix}.
		\end{align*}
		In particular, $\dim(T_0) = 12$.\label{second}
		\item If $\Gamma$ is primitive, and exactly one of $\Gamma$ or $\overline{\Gamma}$ is triangle-free, then  the block dimension decomposition of $T_0$ is 
		\begin{align*}
			\begin{bmatrix}
				1 & 1 & 1\\
				1 & 2 & 2\\
				1 & 2 & 3
			\end{bmatrix}
			\mbox{ or}
			\begin{bmatrix}
				1 & 1 & 1\\
				1 & 3 & 2\\
				1 & 2 & 2
			\end{bmatrix}.
		\end{align*}
		In particular, $\dim(T_0) = 14$.\label{third}
		\item If $\Gamma$ is primitive, and both $\Gamma$ and $\overline{\Gamma}$ admit triangles, then  the block dimension decomposition of $T_0$ is 
		\begin{align*}
			\begin{bmatrix}
				1 & 1 & 1\\
				1 & 3 & 2\\
				1 & 2 & 3
			\end{bmatrix}.
		\end{align*}
		In particular, $\dim(T_0) = 15$.\label{fourth}
	\end{enumerate}
	\label{prop:triply-transitive}
\end{prop}

\subsection{Triply regular graphs}

Let $\Gamma = (\Omega,E)$ be a strongly-regular graph as in Hypothesis~\ref{hyp0}. Recall that $\Gamma$ is {triply regular} if for any $\alpha,\beta,\gamma\in \Omega$, the integer
\begin{align*}
	|\left\{ v\in \Omega: d(\alpha,v) = i,d(\beta,v) = j,d(\gamma,v)=k \right\}|
\end{align*}
is a constant depending only on $i,j,k$ and $a = d(\alpha,\beta), b =  d(\beta,\gamma),c = d(\alpha,\gamma)$. In this case, we let 
\begin{align*}
	\mathrm{p}_{ijk}^{abc} = |\left\{ v\in \Omega: d(\alpha,v) = i,d(\beta,v) = j,d(\gamma,v)=k \right\}|.
\end{align*}
All subconstituents of $\Gamma$ are regular. If $\Gamma$ is triply regular, then its first subconstituent with respect to any vertex, has the property that $\mathrm{p}_{111}^{111}$ is the number of common neighbours of adjacent vertices. Similarly, $\mathrm{p}_{111}^{112}$ is the number of common neighbours of non-adjacent vertices in the first subconsituent. Consequently, the first subconstituent of $\Gamma$ with respect to any vertex is strongly regular. Using the same argument, the second subconstituent of $\Gamma$ with respect to any vertex is also strongly regular.

Let $A_1$ and $A_2$ be respectively the adjacency matrices of $\Gamma$ and $\overline{\Gamma}$. 
Recall that 
\begin{align*}
	T_{0,\omega} = \operatorname{Span}\left\{ E_{i,\omega}^*A_jE_{k,\omega}^*: 0\leq i,j,k\leq 2 \right\}
\end{align*}
for $\omega\in \Omega$.

\begin{thm} \cite{Munemasa}
	The following statements are equivalent.
	\begin{enumerate}[(i)]
		\item The graph $\Gamma$ is triply regular.
		\item There exists $\omega\in \Omega$ such that $T_{0,\omega}  = T_{\omega}$.
		\item For any $\omega\in \Omega$, we have $T_{0,\omega}  = T_{\omega}$.
	\end{enumerate}
	\label{thm:triply-regular}
\end{thm}

In fact, Mumenasa \cite[Proposition~6]{Munemasa} shows an equivalent result to triple regularity that is more practical than Theorem~\ref{thm:triply-regular}.
\begin{thm}\cite[Proposition~6]{Munemasa}
	The graph $\Gamma$ is triply regular if and only if its subconstituents with respect to any vertex are strongly regular.\label{thm:triply-regular-equivalence}
\end{thm}

Another way of characterizing triply-regular strongly-regular graphs is through their Krein parameters. We first recall the following result which can be found in \cite[Theorem~10.7.1]{godsil2013algebraic}.
\begin{lem}
	Let $\Gamma$ be a primitive strongly-regular graph with parameter $(n,k,\lambda,\mu)$. Assume that the eigenvalues of $\Gamma$ are $k>\theta>\tau$. Then,
	\begin{align*}
		\mathrm{q}_{11}^1 &= \theta\tau^2-2\theta^2\tau- \theta^2-k\theta+k\tau^2+2k\tau\geq 0, \mbox{ and }\\
		\mathrm{q}_{22}^2 & =\theta^2\tau - 2\theta\tau^2-\theta^2-k\tau+k\theta^2+2k\theta\geq 0.
	\end{align*}
	If either inequalities is tight, then one of the following is true:
	\begin{enumerate}[(a)]
		\item $\Gamma$ is a $5$-cycle,\label{a} 
		\item Either $\Gamma$ or $\overline{\Gamma}$ has all its first subconstituents empty, and all its second subconstituents strongly regular,\label{b} 
		\item All subconstituents of $\Gamma$ are strongly regular.\label{c} 
	\end{enumerate}
	\label{lem:srg-Krein-parameters}
\end{lem}

Next, we give a partial converse (Lemma~\ref{lem:srg-subconstituent}) to this statement. Before doing this, we need to recall the strongly-regular graphs known as \itbf{Smith graphs}. A strongly regular graph is called a Smith graph if its parameters $(n,k,\lambda,\mu)$ and its eigenvalues $k>\theta>\tau$ satisfy
\begin{align*}
	n &= \frac{2(\theta-\tau)^2 \left((2\theta+1)(\theta-\tau)-3\theta(\theta+1)\right)}{(\theta-\tau)^2-\theta^2(\theta+1)^2}\\
	k &= \frac{-\tau \left((2\theta+1)(\theta-\tau)-\theta(\theta+1)\right)}{(\theta-\tau)+\theta(\theta+1)}\\
	\lambda&= \frac{-\theta(\tau +1)\left((\theta-\tau)-\theta(\theta+3)\right)}{(\theta-\tau)+\theta(\theta+1)}\\
	\mu &= \frac{-(\theta+1)\tau\left((\theta-\tau)-\theta(\theta+1)\right)}{(\theta-\tau)+\theta(\theta+1)}
\end{align*}
where $\theta-\tau\geq \theta(\theta+3)$. The next lemma characterizes strongly-regular graphs with vanishing Krein parameters.
\begin{lem}\cite[Theorem~6.6]{cameron1978strongly}
	Let $\Gamma$ be a graph as in Hypothesis~\ref{hyp0}. We have $\mathrm{q}_{ii}^i = 0$ for some $i\in \{1,2\}$ if and only if $\Gamma$ is either a $5$-cycle, a Smith graph or its complement.\label{lem:Krein-parameters-equivalence}
\end{lem}

The following lemma is a partial converse to Lemma~\ref{lem:srg-Krein-parameters}.
\begin{lem}\cite[Theorem~6.5]{cameron1978strongly}
	If $\Gamma$ is as in Hypothesis~\ref{hyp0}, and both its subconstituents are strongly regular, then one of the following is true.
	\begin{enumerate}[(1)]
		\item $\Gamma$ is a $5$-cycle,
		\item $\Gamma$ has negative Latin square or a pseudo Latin square parameters,
		\item $\Gamma$ is a Smith graph or its complement.
	\end{enumerate}
	\label{lem:srg-subconstituent}
\end{lem}

Now, we combine Lemma~\ref{lem:srg-Krein-parameters} and Lemma~\ref{lem:srg-subconstituent} in the next lemma.
\begin{lem}
	Let $\Gamma = (\Omega,E)$ be a strongly-regular graph as in Hypothesis~\ref{hyp0}. If $\mathrm{q}_{ii}^i>0$ for every $i\in \{1,2\}$, and $\Gamma$ is not a negative Latin square or a pseudo Latin square, then $\Gamma$ is not triply regular. In particular, it is not triply transitive.\label{lem:srg-not-triply-regular}
\end{lem}
\begin{proof}
	Since $\mathrm{q}_{ii}^i>0$ for all $i\in \{1,2\}$, by Lemma~\ref{lem:Krein-parameters-equivalence}, we know that $\Gamma$ is not a $5$-cycle, a Smith graph or its complement. Since $\Gamma$ does not have negative Latin square or pseudo Latin square parameters, using the contrapositive of Lemma~\ref{lem:srg-subconstituent}, we deduce that at least one of the subconstituents of $\Gamma$ is not strongly regular. Hence, $\Gamma$ is not triply regular, and so not triply transitive. This completes the proof. 
\end{proof}

It was observed in \cite[pg~25]{brouwer2022strongly} that if $\Gamma$ is a (primitive) strongly-regular graph with strongly regular subconstituents (including cocliques and cliques) with Latin square or negative Latin square parameters, then the only known examples for $\Gamma$ are: the $n\times n$ grid, graphs with Latin square parameters $(4t^2,t(2t-1),t(t-1),t(t-1))$, graphs with negative Latin square parameters $(r^2(r+3)^2,r^3+3r^2+r,0,r^2+r)$ or $(4t^2,t(2t+1),t(t+1),t(t+1))$, for some integer $t\geq 2$. It was conjectured in \cite[Remark~6.8]{cameron1978strongly} that these are the only examples of strongly-regular graphs with Latin square or negative Latin square parameters admitting strongly-regular subconstituents. As rank $3$ graphs are completely classified, a counterexample of this conjecture is not a rank $3$. Hence, we may state the following for rank $3$ graphs.

\begin{lem}
	Let $\Gamma$ be a strongly-regular graph as in Hypothesis~\ref{hyp0}. Assume that:
	\begin{enumerate}[(a)]
		\item $\Gamma$ is a rank $3$ graph.
		\item $\Gamma$ has Latin square or negative Latin square parameters, and is not a grid graph or does not have parameters $(r^2(r+3)^2,r^3+3r^2+r,0,r^2+r)$.
		\item All subconstituents of $\Gamma$ are strongly regular.
	\end{enumerate}
	Then, there exists $t\geq 2$ such that the parameters of $\Gamma$ or its complement is of the form
	\begin{align*}
		(4t^2,t(2t\pm1),t(t\pm1),t(t\pm1)).
	\end{align*}
	\label{lem:rank3-LS-nLS}	 
\end{lem}

\subsection{Clique extensions}
Let $\Gamma = (\Omega,E)$ be a graph. For any integer $m\geq 1$, the \itbf{$m$-clique extension} of $\Gamma$ is the graph with vertex set $\{1,2,\ldots,m\} \times \Omega$, and two vertices $(i,\omega)$ and $(j,\omega^\prime)$ are adjacent if and only if $\left\{\omega,\omega^\prime\right\} \in E$, or $\omega = \omega^\prime$ and $i\neq j$. The following proposition will be useful for us.
\begin{prop}
	Let $\Gamma$ be a strongly-regular graph with parameters $(n,k,\lambda,\mu)$. For any integer $m\geq 2$, the $m$-clique extension of $\Gamma$ is strongly regular if and only if $\Gamma$ is imprimitive.
\end{prop}
\begin{proof}
	First note that the $m$-clique extension of a clique of size $t$ is a clique of size $tm$. Clearly, if $\Gamma$ is disconnected, then it is a union of cliques of size $k+1$, so the $m$-clique extension is a union of cliques of size $m(k+1)$ which is strongly regular.
	
	Next, assume that the $m$-clique extension of $\Gamma$ is strongly regular. Consider two adjacent vertices $(i,\omega)$ and $(j,\omega^\prime)$. By definition of the $m$-extension, either $\left\{\omega ,\omega^\prime\right\} \in E$ or $\omega = \omega^\prime$ and $i\neq j$. Let us consider these two cases.
	\begin{itemize}
		\item Assume that $\{\omega,\omega^\prime\} \in E$. In this case, the set of common neighbours of the vertices $(i,\omega)$ and $(j,\omega^\prime)$ is the set 
		\begin{align*}
			S = \left\{ (i^\prime,\omega): i\neq i^\prime  \right\} \cup \left\{ (j^\prime,\omega^\prime) : j\neq j^\prime \right\} \cup \left\{ (i^\prime,\alpha): 1\leq i^\prime\leq m, \alpha\in \operatorname{N}_\Gamma(\omega) \cap \operatorname{N}_\Gamma(\omega^\prime) \right\}.
		\end{align*}
		Hence, 
		\begin{align*}
			|S| = 2(m-1)+\lambda m.
		\end{align*}
		\item Assume that $\omega = \omega^\prime$ and $i\neq j$. In this case, the set of common neighbours of the vertices $(i,\omega)$ and $(j,\omega)$ is 
		\begin{align*}
			S^\prime = \left\{ (i^\prime,\omega): i^\prime \neq i, j \right\} \cup \left\{ (i^\prime,\alpha): i^\prime \in \{1,2,\ldots,m\}, \{\alpha,\omega\} \in E \right\}.
		\end{align*}
		Hence, 
		\begin{align*}
			|S^\prime| = m-2+km.
		\end{align*}
	\end{itemize}
	Since $\Gamma$ is strongly regular, we must have $|S| = |S^\prime|$, which implies that $k =\lambda +1$. From what we have seen in \eqref{eq}, the condition $k = \lambda+1$ is a condition for imprimitivity. This completes the proof.
\end{proof}
\begin{cor}
	If $\Gamma = (\Omega,E)$ is a primitive strongly-regular graph, then for any integer $m\geq 2$, the $m$-clique extension of $\Gamma$ is not strongly regular.\label{cor:m-extension}
\end{cor}

\section{The centralizer algebra of a point stabilizer}\label{Sec:TildeT}
Let $\Gamma$ be a strongly-regular graph as in Hypothesis~\ref{hyp0}. Let $G = \Aut{\Gamma}$. Recall that $\tilde{T}_\omega$ is the centralizer algebra of $G_\omega$. That is, the set of all matrices in $\operatorname{M_{|\Omega|}}(\mathbb{C})$ that commute with all permutation matrices $P_{g}$ corresponding to the permutations $g\in G_\omega$. The dimension of $\tilde{T}_\omega$ is well known to be the number of orbitals of $G_\omega$, i.e., the number of orbits of $G_\omega$ on $\Omega \times \Omega$. We will only focus on the case where $G$ is transitive on $\Omega$. Hence, we use $\tilde{T}$ instead of $\tilde{T}_\omega$ for convenience. 

For any subgroup $K\leq G$, define $N_K$ to be the number of orbits of $K$ on $\Omega$.
\begin{lem}
	Assume that $G = \Aut{\Gamma}$ is transitive. Let $\omega_0 = \omega$, $\omega_1\in \Delta_1$, and $\omega_2\in \Delta_2$. Then, 
	\begin{align*}
		\dim(\tilde{T}) = \sum_{i=0}^2 N_{G_\omega \cap G_{\omega_i}}.
	\end{align*}\label{lem:dimension-of-T-tilde}
\end{lem}
\begin{proof}
	We know that $\dim (\tilde{T})$ is the number of orbitals of $G_\omega$ on $\Omega$. Consider the orbital $O = (\alpha,\beta)^{G_\omega}$. By transitivity of $G_\omega$ on $\Delta_i$ for $0\leq i\leq 2$, the orbital of the form $(\alpha,\beta)^{G_\omega}$ is equal to $(\omega,\gamma)^{G_\omega}$, $(\omega_1,\gamma)^{G_\omega}$, or $(\omega_2,\gamma)^{G_\omega}$.
	
	Assume that $O = (\omega,\gamma)^{G_\omega}$, for some $\gamma\in \Omega$. The orbitals of $G_\omega$ are in one-to-one correspondence with the orbits of $G_\omega$. In particular, the orbit $\gamma^{G_\omega}$ corresponds to orbital $(\omega,\gamma)^{G_\omega}$. Hence, there are exactly $N_{G_\omega}$ orbitals of this form ($N_{G_\omega}$ is the rank of $G$).
	
	Assume that $O = (\omega_1,\gamma)^{G_\omega}$, for some $\gamma\in \Omega$. Note that the stabilizer of $\omega_1$ in $G_\omega$ is the group $G_\omega \cap G_{\omega_1}$. If $(\omega_1,\delta)^{G_\omega} = (\omega_1,\delta^\prime)^{G_\omega}$, then there exists $g\in G_\omega$ such that $(\omega_1,\delta)^g =(\omega_1,\delta^\prime)$. In other words, $g\in G_\omega \cap G_{\omega_1}$, and $\delta^g = \delta^\prime$ so $\delta$ and $\delta^\prime$ are in the same orbit of $G_\omega\cap G_{\omega_1}$. Therefore, the number of orbitals of $G_\omega$ of the form $(\omega_1,\delta)^{G_\omega}$, where $\delta\in \Omega$ is equal to the number of distinct orbits of the group $G_\omega\cap G_{\omega_1}$, that is, $N_{G_\omega \cap G_{\omega_1}}$.
	
	Similarly, we use the same argument as above to show that the number of orbitals of the form $O = (\omega_2,\gamma)^{G_\omega}$, for some $\gamma\in \Omega$, is equal to $N_{G_\omega \cap G_{\omega_2}}$.
	This completes the proof. 
\end{proof}

Next, we show that the block dimension decomposition of $\tilde{T}$ is symmetric for rank $3$ graphs.
\begin{thm}
	Let $\Gamma = (\Omega,E)$ be a strongly-regular graph as in Hypothesis~\ref{hyp0}, and assume that $G = \Aut{\Gamma}$ is transitive of rank $3$. Let $\Gamma_1$ and $\Gamma_2$ be its subconstituents with respect to $\omega \in \Omega$.
	Then, the block dimension decomposition of $\tilde{T}$ is of the form
	\begin{align*}
		\begin{bmatrix}
			1 & 1 & 1\\
			1 & r_1 & t\\
			1 & t & r_2
		\end{bmatrix}
	\end{align*}
	for certain integers $r_1,r_2,t\geq 1$. In particular, the integers $r_1$ and $r_2$ are respectively the ranks of the group actions $G_\omega$ on $\Delta_1$ and $G_\omega$ on $\Delta_2$. 
	\label{thm:decomposition-T-tilde}
\end{thm}
\begin{proof}
	Assume $\alpha,\beta \in \Omega$ are such that $\alpha\in \Delta_1$ and $\beta \in \Delta_2$.  Let $r_1$ and $r_2$ be respectively the number of orbitals of $G_\omega$ in $\Delta_1 \times \Delta_1$ and $\Delta_2 \times \Delta_2$. It is clear that $r_1$ and $r_2$ are respectively the $(1,1)$-entry and the $(2,2)$-entry of the block dimension decomposition.
	The number of orbitals of $G_\omega$ in $\Delta_0\times \Delta_i$ is exactly the number of orbits of $G_\omega$ on $\Delta_i$, for $0\leq i\leq 2$. Hence, the first row of the block dimension decomposition of $\tilde{T}$ consists of $1$s since $G$ has rank $3$. Similarly, the number of orbitals in $\Delta_i\times \Delta_0$ is equal to $1$, for any $0\leq i\leq 2$. Consequently, the first column of the block dimension decomposition consists of $1$s.
	
	There is a one-to-one correspondence between the set $\Sigma_1$ of orbitals of $G_\omega$ in $\Delta_1\times \Delta_2$ and the set $\Sigma_2$ of orbitals in $\Delta_2\times \Delta_1$ given by the map $\Sigma_1\to \Sigma_2$ such that  $(\alpha,\gamma)^{G_\omega} \mapsto (\gamma,\alpha)^{G_\omega}.$ From this, we conclude that the matrix of the block dimension decomposition is symmetric. Hence, the block dimension decomposition matrix of $\tilde{T}$ is of the form
	\begin{align*}
		\begin{bmatrix}
			1 & 1 & 1\\
			1 & r_1 & t\\
			1 & t & r_2
		\end{bmatrix}
	\end{align*}
	where $t$ is the number of orbitals of $G_\omega$ in $\Delta_1\times \Delta_2$.
\end{proof}

We will use the following lemma in the proof of the main results.
\begin{lem}
	If $\Gamma = (\Omega,E)$ is a graph whose automorphism group $\Aut{\Gamma}$ is $2$-transitive, then $\Gamma$ is the complete graph or its complement.\label{lem:2-transitive}
\end{lem}
\begin{proof}
	Let $\omega \in \Omega$. If $G = \Aut{\Gamma}$ is $2$-transitive, then $G_\omega$ is transitive on $\Omega \setminus \{\omega\}$. Assume that the vertex $\omega$ has at least one neighbour in $\Gamma$. Since $G_\omega$ is transitive on $\Omega \setminus \{\omega\}$, it is clear that every vertex in $\Omega \setminus \{\omega\}$ is also adjacent to $\omega$. Hence, $\Gamma$ is a complete graph. If $\omega$ is not adjacent to any vertex in $\Omega \setminus \{\omega\}$, then the graph $\Gamma$ is the empty graph.
\end{proof}

\section{Rank $3$ graphs}\label{Sec:rank3-graphs}

Recall that a rank $3$ graph is a strongly-regular graph whose automorphism group is transitive of rank $3$. Such graphs have been classified and the complete classification can be found in \cite[Chapter~11]{brouwer2022strongly}. In particular, it is given in five theorems in \cite[Theorems 11.3.1, 11.3.2, 11.3.3, 11.3.4, 11.4.1]{brouwer2022strongly}. Hence, for the sake of convenience, we will state a readily available classification as presented in \cite{bamberg2023separating}.

Such graphs can be classified into two categories: those belonging to infinite families and the sporadic examples. The infinite families are given in the next theorem.

{\begin{thm}\label{thm:AlmostSimpleAndAffineGraphs}
		Let $\Gamma$ be a rank $3$ graph. If $\Gamma$ is not one of the sporadic examples listed in Table~\ref{tab2}, then it belongs to one of the families below.
		\begin{enumerate}[\rm (1)]
			\item The triangular graph, $T(n)$, for $n\geqslant 4$; \hfill \cite[1.1.7]{brouwer2022strongly}
			\item the Grassmann graph $J_q(n,2)$ for $n\geq 4$; \hfill \cite[\S 3.5.1]{brouwer2022strongly}
			\item the Paley graph, $P(q)$; \hfill \cite[\S 1.1.9]{brouwer2022strongly}
			\item the Peisert graph, $P^*(p^{2t})$; \hfill \cite[\S 7.3.6]{brouwer2022strongly}
			\item the $n \times n$ grid; \hfill \cite[\S 1.1.8]{brouwer2022strongly}
			\item $\nug_m(2)$, for $m> 3$; \hfill \cite[\S 3.1.6]{brouwer2022strongly}
			\item $E_{6,1}(q)$; \hfill \cite[\S 4.9]{brouwer2022strongly}
			\item the collinearity graph of a finite classical polar space or the dual of a finite classical polar space, with rank at least $2$ or rank exactly $2$, respectively; \hfill \cite[Thm 2.2.12 \& 2.2.19]{brouwer2022strongly}
			\item $\vo^\epsilon_{2m}(q)$; \hfill \cite[\S 3.3.1]{brouwer2022strongly} 
			\item the affine half spin graph, $\vd_{5,5}(q)$; \hfill \cite[\S 3.3.3]{brouwer2022strongly}
			\item $\vsz(q)$, for 
			$q=2^{2e-1}$.  \hfill \cite[\S 3.3.1]{brouwer2022strongly} 
			\item the alternating forms graph, $\mathrm{Alt}(5,p^m)$; \hfill \cite[\S 3.4.2]{brouwer2022strongly} 
			\item the Bilinear forms graph $H_q(2, m)$; \hfill \cite[\S 3.4.1]{brouwer2022strongly} 
			\item the van Lint-Schrijver graph, $\vLS(p,e,t)$; \hfill \cite[\S 7.3.1]{brouwer2022strongly}
			\item a connected component of the distance-2 graph of the dual polar graph arising from a polar space of rank $5$ and order $(q,1)$; \hfill \cite[Thm 2.2.20]{brouwer2022strongly} 
		    \item $\no^\epsilon_{2m}(q)$, for $\epsilon = \pm 1$, $m\geq 3$, and $q \in \{2, 3\}$; \hfill \cite[\S 3.1.2]{brouwer2022strongly} 
			\item $\no^{\epsilon}_{2m+1}(q)$, for $\epsilon = \pm 1$, $m\geq 2$ and $q \in \{3, 4, 8\}$. \hfill \cite[\S 3.1.4]{brouwer2022strongly} 
		\end{enumerate}\label{thm:rank3-families}
\end{thm}}

\begin{table}
	\footnotesize 
	\begin{longtable}{cccc} 
		\hline 
		& Graph & Parameters & Reference \\
		\hline \hline
		1 &  $S_{8}$-graph & $(35, 16, 6, 8)$ & \cite[\S10.13]{brouwer2022strongly} \\ 
		2 & $G_2(2)$-graph & $(36,14,4,6)$ &\cite[\S10.14]{brouwer2022strongly} \\ 
		3 & Hoffman-Singleton & $(50, 7, 0, 1)$ & \cite[\S10.19]{brouwer2022strongly} \\ 
		4 & Gewirtz & $(56, 10, 0, 2)$ & \cite[\S10.20]{brouwer2022strongly}\\
		5 & $M_{22}$-graph & $(77, 16, 0, 4)$ &  \cite[\S10.27]{brouwer2022strongly}\\ 
		6 & Higman-Sims & $(100, 22, 0, 6)$ &  \cite[\S10.31]{brouwer2022strongly}\\ 
		7 & Hall-Janko & $(100, 36, 14, 12)$ &  \cite[\S10.32]{brouwer2022strongly}\\
		8 & $S_{10}$-graph  & $(126, 25, 8, 4)$ &  \cite[\S10.40]{brouwer2022strongly} \\
		9 & $U_4(3)$-graph & $(162, 56, 10, 24)$ & \cite[\S10.48]{brouwer2022strongly} \\ 
		10 & Action of $S_4$ on $\rm{PG}(1, 13)$  & $(169, 72, 31, 30)$ &  OA Graph\\ 
		11 & $M_{22}$-graph & $(176, 70, 18, 34)$ &  \cite[\S10.51]{brouwer2022strongly}\\ 
		12 &Berlekamp-van Lint-Seidel & $(243, 22, 1, 2)$ &  \cite[\S10.55]{brouwer2022strongly}\\ 
		13 & Delsarte dual of Berlekamp-van Lint-Seidel & $(243,110,37,60)$ &  \cite[\S10.55]{brouwer2022strongly} \\ 
		14 & $M_{23}$-graph & $(253, 112, 36, 60)$ &  \cite[\S10.56]{brouwer2022strongly} \\ 
		15 & $2^8.S_{10}$-graph & $(256, 45, 16, 6)$ &   \cite[\S10.57]{brouwer2022strongly}\\ 
		16 & $2^8.(A_8 \times S_3)$-graph & $(256, 45, 16, 6)$ &  \cite[\S10.57]{brouwer2022strongly}\\ 
		17 & $2^8.L_2(17)$-graph & $(256,102,38,42)$ &  \cite[\S10.58]{brouwer2022strongly} \\ 
		18 & McLaughlin & $(275, 112, 30, 56)$ &  \cite[\S10.61]{brouwer2022strongly}\\ 
		19 & Action of $S_4$ on $\rm{PG}(1,17)$  & $(289, 96, 35, 30)$ &  OA Graph\\ 
		20 & Action of $S_4$ on $\rm{PG}(1,19)$  & $(361, 144, 59, 56)$ &  OA Graph\\ 
		21 &  $G_2(4)$-graph & $(416, 100, 36, 20)$ &  \cite[\S10.68]{brouwer2022strongly}\\ 
		22 & Sporadic Peisert  & $(529, 264, 131, 132)$ &   \cite[\S10.70]{brouwer2022strongly}\\
		23 & Liebeck 
		& $(625, 144, 43, 30)$ & \cite[\S10.73A]{brouwer2022strongly}\\ 
		24 & Liebeck
		& $(625, 240, 95, 90)$ & \cite[\S10.73B]{brouwer2022strongly}\\ 
		25 & Action of $S_4$ on $\rm{PG}(1,27)$  & $(729, 104, 31, 12)$ &  OA Graph\\ 
		26  & Action of $S_4$ on $\rm{PG}(1,29)$  & $(841, 168, 47, 30)$ &  OA Graph\\ 
		27 & Action of $S_4$ on $\rm{PG}(1,31)$  & $(961, 240, 71, 56)$ & OA Graph \\ 
		28  & Action of $A_5$ on $\rm{PG}(1,31)$ & $(961, 360,139,132)$ &  OA Graph\\
		29 & Dodecad-graph & $(1288, 792, 476, 504)$ & \cite[\S10.80]{brouwer2022strongly}\\ 
		30 & Conway & $(1408, 567, 246, 216)$ &  \cite[\S10.81]{brouwer2022strongly}\\ 
		31 & Action of $A_5$ on $\rm{PG}(1,41)$ & $(1681, 480, 149, 132)$ &   OA Graph\\
		32 & Suzuki & $(1782, 416, 100, 96)$ &  \cite[\S10.83]{brouwer2022strongly} \\ 
		33 & $2^{11}.M_{24}$-graph & $(2048, 276, 44, 36)$ &  \cite[\S10.84]{brouwer2022strongly}\\ 
		34& $2^{11}.M_{24}$-graph & $(2048, 759, 310, 264)$ &  \cite[\S10.85]{brouwer2022strongly}\\ 
		35 & Action of $S_4$ on $\rm{PG}(1,47)$  & $(2209, 1104, 551, 552)$ &  OA Graph \\ 
		36 & Conway   & $(2300, 891, 378, 324)$ &  \cite[\S10.88]{brouwer2022strongly}\\
		37 & $7^4:(6.O_5(3))$-graph 
		& $(2401, 240, 59, 20)$ &  \cite[\S10.89A]{brouwer2022strongly}\\ 
		38 & $7^4:(6.(2^4:S_5))$-graph
		& $(2401, 480, 119, 90)$ &  \cite[\S10.89B]{brouwer2022strongly}\\ 
		39 & $7^4:(3 \times 2.S_7)$-graph
		& $(2401,720,229,210)$ &  \cite[\S10.89C]{brouwer2022strongly}\\
		40 & Action of $A_5$ on ${\rm PG}(1, 49)$
		& $(2401,960,389,380)$ &  \cite[\S10.89D]{brouwer2022strongly}\\ 
		41 & $\fis_{22}$-graph & $(3510, 693, 180, 126)$ &  \cite[\S10.90]{brouwer2022strongly}\\
		42 & Rudvalis & $(4060, 1755, 730, 780)$ &  \cite[\S10.91]{brouwer2022strongly}\\ 
		43 & $2^{12}.HJ.S_3$-graph & $(4096, 1575, 614, 600)$ &  \cite[\S10.92]{brouwer2022strongly} \\ 
		44 & Action of $A_5$ on $\rm{PG}(1,71)$ & $(5041, 840, 179, 132)$&  OA Graph\\
		45 & Action of $A_5$ on $\rm{PG}(1,79)$ & $(6241, 1560, 419, 380)$ &  OA Graph\\
		46 & $3^8.2^{1+6}.O_6^-(2).2$-graph & $(6561, 1440, 351, 306)$ &  \cite[\S10.93]{brouwer2022strongly}\\ 
		47 & Action of $A_5$ on $\rm{PG}(1,89)$& $(7921, 2640, 899, 870)$ &  OA Graph\\
		48 & $\fis_{22}$-graph & $(14080, 3159, 918, 648)$ &  \cite[\S10.94]{brouwer2022strongly}\\ 
		49 & $5^6.4.\mathit{HJ}.2$-graph & $(15625, 7560, 3655, 3660)$ &  \cite[\S10.95]{brouwer2022strongly}\\ 
		50 & $\fis_{23}$-graph & $(31671, 3510, 693, 351)$ &  \cite[\S10.96]{brouwer2022strongly}\\
		51 & $\fis_{23}$-graph & $(137632,28431, 6030, 5832)$ &  \cite[\S10.97]{brouwer2022strongly} \\ 
		52 & $\fis_{24}$-graph & $(306936, 31671, 3510, 3240)$ &  \cite[\S10.99]{brouwer2022strongly}\\ 
		53 & Suz-graph & $(531441,65520,8559,8010)$ &  \cite[\S10.100]{brouwer2022strongly}\\
		\hline \caption{Sporadic rank $3$ graphs}\label{tab2}
	\end{longtable}
\end{table}

\section{Triply-transitive graphs}\label{Sec:first-cases}

Let $\Gamma = (\Omega,E)$ be a strongly-regular graph as in Hypothesis~\ref{hyp0}. Let $T = T_\omega$ be the Terwilliger algebra of $\Gamma$ with respect to a vertex $\omega\in \Omega$, and let $T_0 = T_{0,\omega}$ be the subspace defined in \eqref{eq:T0}. Let $\tilde{T} = \tilde{T}_\omega$ be the centralizer algebra of $G_\omega$ with respect to this vertex $\omega\in \Omega$. In this section, we consider the cases where $\Gamma$ is triply transitive (see Definition~\ref{def:1}).

First, we will determine some properties of the automorphism groups of such graphs in Section~\ref{subsect:automorphism}. In particular, we show that a triply-transitive graph must be a rank $3$ graph. Next, in order to classify triply-transitive graphs, we use the fact that they must be triply regular. In particular, we will use Proposition~\ref{prop:triply-transitive} to analyze all possible cases by considering the dimensions of $T_0$. The cases \eqref{first},\eqref{second},\eqref{third}, and \eqref{fourth} in Proposition~\ref{prop:triply-transitive} are respectively considered in Sections~\ref{subsect:Imprimitive}, \ref{subsect:five-cycles}, \ref{subsect:triangle-free}, and \ref{subsect:general}.

\subsection{Automorphism groups of triply-transitive graphs}\label{subsect:automorphism}
We prove that if $\Gamma$ is triply transitivity, then it must be a rank $3$ graph.  
{In the next theorem, we prove a result that is slightly more general than what is needed for triply transitive graphs. To be precise, we determine the properties of the automorphism group of graphs satisfying the relation $T_0 =T = \tilde{T}$, with respect to $\omega\in \Omega$.}
\begin{thm}
	Let $\Gamma=(\Omega,E)$ be a strongly-regular graph as in Hypothesis~\ref{hyp0}. If $\Gamma = (\Omega,E)$ is such that $T_0 = T = \tilde{T}$ with respect to a vertex $\omega\in \Omega$, then the automorphism group $G=\Aut{\Gamma}$ satisfies one of the following.
	\begin{enumerate}[(a)]
		\item  $G$ is transitive of rank $3$ on $\Omega$.\label{first-transitive}
		\item  $G = G_\omega$, and $G$ is transitive on $\Delta_0 = \{\omega\}$, $\Delta_1$ and $\Delta_2$.\label{second-intransitive}
	\end{enumerate}
	\label{thm:triple-transitive-rank3}
\end{thm}
\begin{proof}
	Since $T_0 = T = \tilde{T}$, we know that the first row of the block dimension decomposition of $\tilde{T}$ consists of $1$s. This implies that $G_\omega$ has three orbits on $\Omega$, which are $\Delta_0 = \{\omega\}$, the neighbourhood $\Delta_1$ of $\omega$ in $\Gamma$, and the set of non-neighbours $\Delta_2$ of $\omega$. Hence, if $G$ is transitive, then it must have rank $3$, which proves \eqref{first-transitive}. 
	
	Now, assume that $G$ is intransitive on $\Omega$ and let $\Omega_1,\Omega_2,\ldots,\Omega_d$ be the orbits of $G$ on $\Omega$, for some $d\geq 2$. If $d\geq 4$, then there exist two vertices $\alpha$ and $\beta$ which are both in $\Delta_1$ or both in $\Delta_2$ such that $\alpha \in \Omega_i$ and $\beta \in \Omega_j$ for some distinct $1\leq i,j\leq d$. Since $G_\omega$ is transitive on $\Delta_1$ and $\Delta_2$, we conclude that there exists an element of $G_\omega \leq G$ which maps $\alpha$ to $\beta$. This is a contradiction to the fact that $\Omega_i$ and $\Omega_j$ are distinct orbits of $G$ containing $\alpha$ and $\beta$, respectively. Consequently, $d\leq 3$.
	
	Next, we show that the case $d=2$ is not possible. Assume by contradiction that $d = 2$. Then, it is not hard to see that one of the following cases may occur.
	\begin{enumerate}
		\item $\Omega_1 = \Delta_0$ and $\Omega_2 = \Delta_1\cup \Delta_2$,
		\item $\Omega_1 = \Delta_0 \cup \Delta_2$ and $\Omega_2 = \Delta_1$,
		\item $\Omega_1 = \Delta_0 \cup \Delta_1$ and $\Omega_2 = \Delta_2$.
	\end{enumerate}
	If (i) holds, then $G = G_\omega$. If $\alpha \in \Delta_1$ and $\beta \in \Delta_2$, then using the fact that $\Omega_2 = \Delta_1\cup \Delta_2$ is an orbit of $G$, there exists $g\in G$ such that $\alpha^g = \beta$. Since $\{\omega,\alpha\} \in E$ and $\{ \omega,\beta \} \not\in E$, we have $\{\omega,\alpha\}^g = \{ \omega,\beta \} \in E$, thus contradicting the fact that $g \in G$ is an automorphism. Hence, (i) is impossible.
	
	If (ii) holds, then choose $\beta \in \Delta_2$. Since $\omega,\beta \in \Omega_1$, there exists $g\in G$ such that $\omega^g = \beta$. Since $\omega $ is adjacent to every vertex in $ \Delta_1$, we conclude that $\omega^g = \beta $ is also adjacent to $ \Delta_1^g = \Delta_1$. Consequently, every vertex in $\Delta_2$ is adjacent to all vertices of $\Delta_1$, which implies that $\Gamma$ is imprimitive. This is a contradiction to the fact that $G$ is intransitive since the automorphism group of an imprimitive strongly-regular graph is transitive. This rules out case (ii). 
	
	If (iii) holds, then we consider the complement of $\Gamma$. Note that $\Gamma$ and its complement have the same automorphism group. For the complement, $\Omega_1$ is the set of vertices not adjacent to $\omega$, and so the orbit $\Omega_1$ is the union of $\Delta_0 = \{\omega\}$ and $\Delta_1$ consisting of the non-neighbours. Therefore, applying the same argument as in (ii), we conclude that case (iii) is not possible.
	
	Overall, we have shown that $d$ cannot be equal to $2$. Therefore, $d=3$, and without loss of generality, we can assume that $\Omega_1=\Delta_0 = \{\omega\}$, $\Omega_2=\Delta_1$, and $\Omega_3 = \Delta_2$.
\end{proof}

We do not have any example of strongly-regular graphs for which \eqref{second-intransitive} in Theorem~\ref{thm:triple-transitive-rank3} holds. Therefore, we ask the following question.

\begin{qst}
	Is there a strongly-regular graph $\Gamma$ for which Theorem~\ref{thm:triple-transitive-rank3}~\eqref{second-intransitive} holds?
\end{qst}
\begin{rmk}
	We note that there are examples of graphs $\Gamma$ satisfying Hypothesis~\ref{hyp0} for which there exists $\omega \in \Omega$ such that $G = \Aut{\Gamma} = G_\omega$, and $G$ is transitive on $\Delta_0 = \{\omega\}$, $\Delta_1$ and $\Delta_2$. For instance, the association scheme \verb*|as25n10| (see \cite{Hanaki}) is a non-Schurian $2$-class scheme with automorphism group $G$ isomorphic to \verb|(C3 x A4):C2|. The group $G$ is intransitive, and its orbits are
	\begin{align*}
		\Delta_0 &= \{4\},\\ \Delta_1 &= \{1, 2, 8, 3, 15, 9, 14, 22, 16, 10, 20, 21\},\\\Delta_2 &= \{5, 6, 13, 7, 18, 12, 19, 23, 17, 11, 24, 25\}.
	\end{align*}
	However, we do not have $T_0 = T = \tilde{T}$ with respect to vertex $4$, since $\dim(T_0) = 15$ and $\dim(T) = \dim(\tilde{T}) = 19 $.
\end{rmk}

	The following corollary immediately follows from Theorem~\ref{thm:triple-transitive-rank3}.
	\begin{cor}
		If $\Gamma$ is triply transitive, then it is a rank $3$ graph.\label{cor:rank3}
	\end{cor}

Therefore, we may define all related subspaces of $\Gamma$ with respect to a fixed vertex $\omega\in \Omega$, whenever $\Gamma$ is triply transitive. For convenience, we consider the following assumption.

\begin{hyp}
	Let $\Gamma = (\Omega,E)$ be a strongly-regular graph as in Hypothesis~\ref{hyp0}. Assume further that $\Gamma$ is a rank $3$ graph with full automorphism group $G = \Aut{\Gamma}$. Fix $\omega\in \Omega$ and consider the Terwilliger algebra $T$ with respect to $\omega$, the centralizer algebra $\tilde{T}$ of the point stabilizer $G_\omega$, and the subspace $T_0$ with respect to $\omega$. Suppose that $\Gamma$ is triply transitive.
	\label{hyp}
\end{hyp}

\subsection{Imprimitive strongly-regular graphs}\label{subsect:Imprimitive}
In this section, we determine the strongly-regular graphs that are both triply transitive and imprimitive. Assume that $\Gamma$ is imprimitive and triply transitive. By imprimitivity, we know that $\Gamma$ is a complete $k$-partite graph, for some $k\geq 2$. By Proposition~\ref{prop:triply-transitive}~\eqref{first}, we know that $\dim(T_0)  \in \{11,12\}$ for such graphs depending on $k$. 
In order to determine the triply transitive and imprimitive strongly-regular graphs, it is enough to determine the values of $k$ for which $\dim(\tilde{T}) = \dim(T_0)$.

\begin{lem}\label{imprimitive Ttilde}
	Assume that ${\Gamma}$ is a complete $k$-partite graph. Then, we have 
	\begin{align*}
		\dim (\tilde{T}) = 
		\begin{cases}
			11 & \mbox{ if  }k =2,\\
			12 & \mbox{ otherwise.}\\
		\end{cases}
	\end{align*} 
	In particular, $\Gamma$ is triply transitive for all $k\geq 2$.
\end{lem}
\begin{proof}
	Let $B_1,\ldots,B_k$ be the parts of $\Gamma$, where $|B_i| = n$ for all $1\leq i\leq k$. Assume that $\omega \in B_1$. Then clearly, $\Delta_1 = B_2\cup B_3\ldots \cup B_k$ and $\Delta_2 = B_1\setminus\{\omega\}$.
	
	 The full automorphism of $\Gamma$ is $G = \sym{n}\wr \sym{k}$. We have $G_\omega = \sym{n-1}\times \left(\sym{n}\wr\sym{k-1}\right)$. Let us compute the number of orbitals of $G_\omega$ in its action on $\Omega$.
	
	Assume that $k\geq 3$. Recall that the orbits of $G_\omega$ are $\Delta_0 = \{\omega\}, \Delta_1$, and $\Delta_2.$ The number of orbitals of $G_\omega$ is equal to the number of orbits of $G_\omega\cap G_\alpha$, where $\alpha$ runs through every representative of $\Delta_0, \Delta_1,$ and $\Delta_2$.
	\begin{itemize}
		\item Since $|\Delta_0| = 1$, there are $3$ orbits for $G_\omega$.
		\item Given $\alpha\in \Delta_1$, the subgroup $G_\omega\cap G_\alpha$ has one orbit on $\Delta_0$, two orbits on $\Delta_1$, and one orbit on $\Delta_2$.
		\item Given $\alpha\in \Delta_2$, the subgroup of $G_\omega\cap G_\alpha$ has one orbit on $\Delta_0$, one orbit on $\Delta_2,$ and three orbits on $\Delta_2$.
	\end{itemize}
	Consequently, we have $\dim (\tilde{T}) = 12$.
	
	Next, we consider the case where $k = 2$ or equivalently $\Gamma$ is a complete bipartite graph. In this case, we have the following.
	\begin{itemize}
		\item Since $|\Delta_0| = 1$, there are $3$ orbits for $G_\omega$.
		\item Given $\alpha\in \Delta_1$, the subgroup $G_\omega\cap G_\alpha$ has one orbit on $\Delta_0$, two orbits on $\Delta_1$, and one orbit on $\Delta_2$.
		\item Given $\alpha\in \Delta_2$, the subgroup of $G_\omega\cap G_\alpha$ has one orbit on $\Delta_0$, one orbit on $\Delta_2,$ and two orbits on $\Delta_2$.
	\end{itemize}
	Consequently, we have $\dim(\tilde{T}) = 11$.
	
	The second statement follows immediately from Proposition~\ref{prop:triply-transitive}\eqref{first}.
\end{proof}

\begin{thm}
	All imprimitive strongly-regular graphs are triply transitive.\label{thm:imprimitive}
\end{thm}

\subsection{$\Gamma$ and its complement are triangle free}\label{subsect:five-cycles}
Now, we determine the triply-transitive graphs with the property that the graph and its complement are both triangle free.
\begin{lem}
	Let $\Gamma$ be a primitive strongly-regular graph as in Hypothesis~\ref{hyp}. If $\Gamma$ and its complement are both triangle free, then $\Gamma$ is the cycle on 5 vertices.
\end{lem}
\begin{proof}
	By Proposition~\ref{prop:triply-transitive}\eqref{second}, we know that $\dim(T_0) = 13$, and the block dimension decomposition of $T_0 = \tilde{T}$ is
	\begin{align*}
		\begin{bmatrix}
			1 & 1 & 1\\
			1 & 2 & 2\\
			1 & 2 & 2
		\end{bmatrix}.
	\end{align*}
	By Theorem~\ref{thm:decomposition-T-tilde}, we know that the rank of $G_\omega$ on $\Delta_1$ and $\Delta_2$ are both equal to $2$. By Lemma~\ref{lem:2-transitive}, the group $G_\omega$ acts $2$-transitively on the first subconstituent of $\Gamma$ (with respect to $\omega$). Similarly, it also acts $2$-transitively on the second subconstituent. We conclude that the first and second subconstituents are cliques or cocliques. Whenever $\Delta_1$ and $\Delta_2$ has more than $2$ elements, then it is clear that one of $\Gamma$ and $\overline{\Gamma}$ has a triangle. Hence, both subconstituents of $\Gamma$ must have size at most $2$.
	
	If $\Gamma$ is a $4$-cycle, then $\Gamma$ is imprimitive. Hence, $\Gamma$ must have five vertices. It is an easy exercise to show that $\Gamma$ must be a cycle on five vertices.
\end{proof}

\subsection{Exactly one of $\Gamma$ or its complement is triangle free}\label{subsect:triangle-free}
Let $\Gamma$ be a strongly-regular graph as in Hypothesis~\ref{hyp}. Assume that $\Gamma$ is triangle-free and its complement admits triangles. By Proposition~\ref{prop:triply-transitive}\eqref{third}, we have $\dim(T_0) = 14$. We may assume without loss of generality that the block dimension decomposition of $\tilde{T}$ is 
\begin{align*}
	\begin{bmatrix}
		1 & 1 & 1\\
		1 & 2 & 2\\
		1 & 2 & 3
	\end{bmatrix}. 
\end{align*}
If $\Gamma$ is imprimitive, then it must be a complete bipartite graph, which has already been considered in Section~\ref{subsect:Imprimitive}. Now, assume that $\Gamma$ is primitive. Using the classification of primitive rank $3$ graphs in Theorem~\ref{thm:rank3-families}, the only triangle-free primitive rank $3$ graphs whose complements admit triangles are:
\begin{itemize}
	\item the Petersen graph (isomorphic to the complement of the Johnson graph $J(5,2)$),
	\item the complement of Clebsch (which is isomorphic to the halved 5-cube),
	\item the Hoffman-Singleton graph,
	\item the Gewirtz graph,
	\item the $\mathieu{22}$-graph,
	\item the Higman-Sims graph.
\end{itemize}
Using \verb*|Sagemath| \cite{sagemath}, we have computed the dimensions of the relevant subspaces, and listed them in the next table.

\begin{table}[H]
	\centering\tiny
	\begin{longtable}{|c|c|c|c|c|c|c|} \hline
		SRG & Name & Automorphism group & $\dim (T_0)$ & $\dim(T)$ & $\dim(\tilde{T})$ & Decomposition \\ \hline
		(10,3,0,1) & Petersen graph & \verb|S5| & $14$ & $15$ & $15$ & $\left(\begin{array}{rrr}
			1 & 1 & 1 \\
			1 & 2 & 2 \\
			1 & 2 & 4
		\end{array}\right)$ \\ \hline
		(16,5,0,2) & Complement of the Clebsch graph & \verb|(C2 x C2 x C2 x C2) : S5| & $14$ & $14$ & $14$ & $\left(\begin{array}{rrr}
			1 & 1 & 1 \\
			1 & 2 & 2 \\
			1 & 2 & 3
		\end{array}\right)$ \\ \hline
		(50,7,0,1) & Hoffman-Singleton graph & \verb|PSU(3,5) : C2| & $14$ & $15$ & $15$ & $\left(\begin{array}{rrr}
			1 & 1 & 1 \\
			1 & 2 & 2 \\
			1 & 2 & 4
		\end{array}\right)$ \\ \hline
		(56,10,0,2) & Gerwirtz graph & \verb|PSL(3,4) : (C2 x C2)| & $14$ & $15$ & $16$ & $\left(\begin{array}{rrr}
			1 & 1 & 1 \\
			1 & 5 & 2 \\
			1 & 2 & 2
		\end{array}\right)$ \\ \hline
		(77,16,0,4) & $\mathieu{22}$-graph & \verb|M22 : C2| & $14$ & $15$ & $16$ & $\left(\begin{array}{rrr}
			1 & 1 & 1 \\
			1 & 5 & 2 \\
			1 & 2 & 2
		\end{array}\right)$ \\ \hline
		(100,22,0,6) & Higman-Sims graph & \verb|HS: C2| & $14$ & $14$ & $14$ & $\left(\begin{array}{rrr}
			1 & 1 & 1 \\
			1 & 2 & 2 \\
			1 & 2 & 3
		\end{array}\right)$ \\ \hline
	\end{longtable}
\end{table}

Therefore, the only triangle-free primitive rank $3$ graphs that are triply transitive are the complement of the Clebsch graph and the Higman-Sims graph. We state this result in the next lemma.
\begin{lem}
	If $\Gamma$ is a primitive strongly-regular graph as in Hypothesis~\ref{hyp} with the property that exactly one of $\Gamma$ or its complement is triangle free, then $\Gamma$ is one of:
	\begin{enumerate}[(i)]
		\item the complement of the Clebsch graph,
		\item the Higman-Sims graph.
	\end{enumerate}
\end{lem}

\section{$\Gamma$ and its complement admit triangles}\label{subsect:general}
In this section, we study the remaining cases that were not considered in the previous section. Therefore, we assume that $\Gamma = (\Omega,E)$ is a strongly-regular graph as in Hypothesis~\ref{hyp}. Moreover, assume that $\Gamma$ and its complement admit triangles. By Proposition~\ref{prop:triply-transitive}\eqref{fourth}, we have $\dim(T_0) = 15$ in this case.

We will use the classification of rank $3$ graphs in Theorem~\ref{thm:rank3-families} to determine the possibilities for the graph $\Gamma$.

\subsection{The sporadic cases of rank $3$ graphs}

First, we deal with most of the rank $3$ graphs that are not part of an infinite family. These are given in Table~\ref{tab2}. The main tools that we use are \verb*|Sagemath| and Lemma~\ref{lem:srg-not-triply-regular}. We use \verb*|Sagemath| for the rank $3$ graphs whose automorphism groups are available in the library of \verb*|PrimitiveGroups|. For the graphs with many vertices, we can sometimes use Lemma~\ref{lem:srg-not-triply-regular}, unless they have Latin square or negative Latin square parameters.
The result of our computations are given in Table~\ref{tab}. In particular, we determine all triply-transitive graphs in Table~\ref{tab}, except for those in Row~43,49,53. \newpage

{\tiny
	\centering
	
	\begin{longtable}{cccccc} 
		\hline 
		& Graph & Parameters & Reference &Triply transitive&Comments\\
		\hline \hline
		1 &  $S_{8}$-graph & $(35, 16, 6, 8)$ & \cite[\S10.13]{brouwer2022strongly} & No& \verb*|Sagemath|\\ 
		2 & $G_2(2)$-graph & $(36,14,4,6)$ &\cite[\S10.14]{brouwer2022strongly} & No&\verb*|Sagemath|\\ 
		3 & Hoffman-Singleton & $(50, 7, 0, 1)$ & \cite[\S10.19]{brouwer2022strongly} & No&\verb*|Sagemath|\\ 
		4 & Gewirtz & $(56, 10, 0, 2)$ & \cite[\S10.20]{brouwer2022strongly}&No&\verb*|Sagemath|\\
		5 & $M_{22}$-graph & $(77, 16, 0, 4)$ &  \cite[\S10.27]{brouwer2022strongly}&No&\verb*|Sagemath|\\ 
		6 & Higman-Sims & $(100, 22, 0, 6)$ &  \cite[\S10.31]{brouwer2022strongly}& Yes&\verb*|Sagemath|\\ 
		7 & Hall-Janko & $(100, 36, 14, 12)$ &  \cite[\S10.32]{brouwer2022strongly}& No&\verb*|Sagemath|\\
		8 & $S_{10}$-graph  & $(126, 25, 8, 4)$ &  \cite[\S10.40]{brouwer2022strongly} &No &\verb*|Sagemath|\\
		9 & $U_4(3)$-graph & $(162, 56, 10, 24)$ & \cite[\S10.48]{brouwer2022strongly} & No&\verb*|Sagemath|\\ 
		10 & Action of $S_4$ on $\rm{PG}(1, 13)$  & $(169, 72, 31, 30)$ &  OA Graph&No&\verb*|Sagemath|\\ 
		11 & $M_{22}$-graph & $(176, 70, 18, 34)$ &  \cite[\S10.51]{brouwer2022strongly}& No &\verb*|Sagemath|\\ 
		12 &Berlekamp-van Lint-Seidel & $(243, 22, 1, 2)$ &  \cite[\S10.55]{brouwer2022strongly}&No &\verb*|Sagemath|\\ 
		13 & Delsarte dual of Berlekamp-van Lint-Seidel & $(243,110,37,60)$ &  \cite[\S10.55]{brouwer2022strongly} &No&\verb*|Sagemath|\\ 
		14 & $M_{23}$-graph & $(253, 112, 36, 60)$ &  \cite[\S10.56]{brouwer2022strongly} & No&\verb*|Sagemath|\\ 
		15 & $2^8.S_{10}$-graph & $(256, 45, 16, 6)$ &   \cite[\S10.57]{brouwer2022strongly}&No&\verb*|Sagemath|\\ 
		16 & $2^8.(A_8 \times S_3)$-graph & $(256, 45, 16, 6)$ &  \cite[\S10.57]{brouwer2022strongly}& No&\verb*|Sagemath|\\ 
		17 & $2^8.L_2(17)$-graph & $(256,102,38,42)$ &  \cite[\S10.58]{brouwer2022strongly} &No &\verb*|Sagemath|\\ 
		18 & McLaughlin & $(275, 112, 30, 56)$ &  \cite[\S10.61]{brouwer2022strongly}&Yes&\verb*|Sagemath|\\ 
		19 & Action of $S_4$ on $\rm{PG}(1,17)$  & $(289, 96, 35, 30)$ &  OA Graph& No &\verb*|Sagemath|\\ 
		20 & Action of $S_4$ on $\rm{PG}(1,19)$  & $(361, 144, 59, 56)$ &  OA Graph& No &\verb*|Sagemath|\\ 
		21 &  $G_2(4)$-graph & $(416, 100, 36, 20)$ &  \cite[\S10.68]{brouwer2022strongly}& No&\verb*|Sagemath|\\ 
		22 & Sporadic Peisert  & $(529, 264, 131, 132)$ &   \cite[\S10.70]{brouwer2022strongly}&No&\verb*|Sagemath|\\
		23 & Liebeck 
		& $(625, 144, 43, 30)$ & \cite[\S10.73A]{brouwer2022strongly}&No &\verb*|Sagemath|\\ 
		24 & Liebeck
		& $(625, 240, 95, 90)$ & \cite[\S10.73B]{brouwer2022strongly}&No&\verb*|Sagemath|\\ 
		25 & Action of $S_4$ on $\rm{PG}(1,27)$  & $(729, 104, 31, 12)$ &  OA Graph&No&\verb*|Sagemath|\\ 
		26  & Action of $S_4$ on $\rm{PG}(1,29)$  & $(841, 168, 47, 30)$ &  OA Graph&No&\verb*|Sagemath|\\ 
		27 & Action of $S_4$ on $\rm{PG}(1,31)$  & $(961, 240, 71, 56)$ & OA Graph &No&\verb*|Sagemath|\\ 
		28  & Action of $A_5$ on $\rm{PG}(1,31)$ & $(961, 360,139,132)$ &  OA Graph&No&\verb*|Sagemath|\\
		29 & Dodecad-graph & $(1288, 792, 476, 504)$ & \cite[\S10.80]{brouwer2022strongly}&No&\verb*|Sagemath|\\ 
		30 & Conway & $(1408, 567, 246, 216)$ &  \cite[\S10.81]{brouwer2022strongly}&No&\verb*|Sagemath|\\ 
		31 & Action of $A_5$ on $\rm{PG}(1,41)$ & $(1681, 480, 149, 132)$ &   OA Graph& No &\verb*|Sagemath|\\
		32 & Suzuki & $(1782, 416, 100, 96)$ &  \cite[\S10.83]{brouwer2022strongly} & No &\verb*|Sagemath|\\ 
		33 & $2^{11}.M_{24}$-graph & $(2048, 276, 44, 36)$ &  \cite[\S10.84]{brouwer2022strongly}& No&\verb*|Sagemath|\\ 
		34& $2^{11}.M_{24}$-graph & $(2048, 759, 310, 264)$ &  \cite[\S10.85]{brouwer2022strongly}& No&\verb*|Sagemath|\\ 
		35 & Action of $S_4$ on $\rm{PG}(1,47)$  & $(2209, 1104, 551, 552)$ &  OA Graph &No&\verb*|Sagemath|\\ 
		36 & Conway   & $(2300, 891, 378, 324)$ &  \cite[\S10.88]{brouwer2022strongly}& No&\verb*|Sagemath|\\
		37 & $7^4:(6.O_5(3))$-graph \footnote{The full automorphism group is the primitive group \sf{PrimitiveGroup}(2401,1049).}
		& $(2401, 240, 59, 20)$ &  \cite[\S10.89A]{brouwer2022strongly}& No & $\dim(\tilde{T}) = 143$\\ 
		38 & $7^4:(6.(2^4:S_5))$-graph
		& $(2401, 480, 119, 90)$ &  \cite[\S10.89B]{brouwer2022strongly}&No&$\dim(\tilde{T}) = 563$\\ 
		39 & $7^4:(3 \times 2.S_7)$-graph
		& $(2401,720,229,210)$ &  \cite[\S10.89C]{brouwer2022strongly}& No& $\dim(\tilde{T}) = 249$\\
		40 & Action of $A_5$ on ${\rm PG}(1, 49)$
		& $(2401,960,389,380)$ &  \cite[\S10.89D]{brouwer2022strongly}&No &Latin square\\ 
		41 & $\fis_{22}$-graph & $(3510, 693, 180, 126)$ &  \cite[\S10.90]{brouwer2022strongly}&No&Lemma~\ref{lem:srg-not-triply-regular}\\
		42 & Rudvalis & $(4060, 1755, 730, 780)$ &  \cite[\S10.91]{brouwer2022strongly}&No&local graph not strongly-regular\\ 
		43 & $2^{12}.HJ.S_3$-graph & $(4096, 1575, 614, 600)$ &  \cite[\S10.92]{brouwer2022strongly} &No&Lemma~\ref{lem:rank3-LS-nLS}\\ 
		44 & Action of $A_5$ on $\rm{PG}(1,71)$ & $(5041, 840, 179, 132)$&  OA Graph&No&Latin square\\
		45 & Action of $A_5$ on $\rm{PG}(1,79)$ & $(6241, 1560, 419, 380)$ &  OA Graph&No&Latin square\\
		46 & $3^8.2^{1+6}.O_6^-(2).2$-graph & $(6561, 1440, 351, 306)$ &  \cite[\S10.93]{brouwer2022strongly}& No& Lemma~\ref{lem:srg-not-triply-regular}\\ 
		47 & Action of $A_5$ on $\rm{PG}(1,89)$& $(7921, 2640, 899, 870)$ &  OA Graph&No&Latin square\\
		48 & $\fis_{22}$-graph & $(14080, 3159, 918, 648)$ &  \cite[\S10.94]{brouwer2022strongly}& No&Lemma~\ref{lem:srg-not-triply-regular}\\ 
		49 & $5^6.4.\mathit{HJ}.2$-graph & $(15625, 7560, 3655, 3660)$ &  \cite[\S10.95]{brouwer2022strongly}& No &Lemma~\ref{lem:rank3-LS-nLS}\\ 
		50 & $\fis_{23}$-graph & $(31671, 3510, 693, 351)$ &  \cite[\S10.96]{brouwer2022strongly}& No  & Lemma~\ref{lem:srg-not-triply-regular} \\
		51 & $\fis_{23}$-graph & $(137632,28431, 6030, 5832)$ &  \cite[\S10.97]{brouwer2022strongly} &No&Lemma~\ref{lem:srg-not-triply-regular}\\ 
		52 & $\fis_{24}$-graph & $(306936, 31671, 3510, 3240)$ &  \cite[\S10.99]{brouwer2022strongly}& No&Lemma~\ref{lem:srg-not-triply-regular}\\ 
		53 & Suz-graph & $(531441,65520,8559,8010)$ &  \cite[\S10.100]{brouwer2022strongly}&No&Lemma~\ref{lem:rank3-LS-nLS}\\
		\hline \caption{Sporadic rank $3$ graphs}\label{tab}
	\end{longtable}
	
}

\subsection{Triangular graphs}
The graph $\Gamma$ is isomorphic to the Johnson graph $J(\Omega,2)$ in this case. Assume that $\Omega = \{1,2,\ldots,n\}$. Its parameters are
\begin{align*}
	\left(\binom{n}{2},2(n-2),n-2,4\right). 
\end{align*}

The first subconstituent of $J(\Omega,2)$ with respect to the vertex $A = \{1,2\}$ is the subgraph of $J(\Omega,2)$ induced by
\begin{align*}
	V = \left\{ \{1,i\}:3\leq i\leq n \right\}\cup \left\{ \{2,i\}: 3\leq i\leq n \right\}.
\end{align*}
The pair of adjacent vertices $(\{1,3\},\{1,4\})$ has $n-4$ common neighbours in $\Delta_1$, whereas the pair of adjacent vertices $(\{1,3\},\{2,3\})$ has no common neighbours in $\Delta_1$. Consequently, the triangular graph is not triply regular by Theorem~\ref{thm:triply-regular-equivalence}, and thus cannot be triply transitive.
\begin{thm}
	For any set $\Omega$ with at least $5$ elements, the Johnson graph $J(\Omega,2)$ is not triply transitive.\label{thm:triangular}
\end{thm}

\subsection{Grassmann graphs}
Let $q$ be a power of a prime and $n\geq d\geq 1$ be two integers. Recall that the \itbf{Grassmann graph} $J_q(n,d)$ is the graph whose vertex set is the set of all $d$-dimensional subspaces of $\mathbb{F}_q^n$, and two vertices $U$ and $V$ are adjacent if $U\cap V$ is a $(d-1)$-dimensional space. 
When $d = 2$, then $J_q(n,2)$ is a strongly-regular graph with parameters
\begin{align}
	\left(\qbinom{n}{2}{q},(q+1)\left(\qbinom{n-1}{1}{q}-1\right),\qbinom{n-1}{1}{q}+q^2-2,(q+1)^2\right),\label{eq:grassman-parameter}
\end{align}
where $\qbinom{n}{d}{q} = \frac{(q^n-1)\ldots(q^{n-d+1}-1)}{(q^d-1)\ldots(q-1)}$. Equivalently, the graph $J_q(n,2)$ is the graph whose vertices consist of the lines of $\pg{n-1}{q}$, and two lines are adjacent if they intersect.

For any $1\leq i\leq n$, we let $e_i$ be the canonical vector of $\mathbb{F}_q^n$.
Let $V = \langle e_1,e_2\rangle, $ and consider the partition $\Delta_0 = \{V\},$ $\Delta_1$ the neighbourhood of $V$ in $J_q(n,2)$, and $\Delta_2$ the non-neighbours. Let $U = \langle e_1,e_3 \rangle, W = \langle e_1,e_4 \rangle, U^\prime = \langle e_1,e_3\rangle$, and $W^\prime = \langle e_2,e_3 \rangle$. It is clear that in the graph $J_q(n,2)$, we have $U,W,U^\prime,W^\prime \in \Delta_1$. Moreover, we have $U\sim W$ and $U^\prime \sim W^\prime$ in the first subconstituent. Let us show that the number of common neighbours of $(U,W)$ and $(U^\prime,W^\prime)$ in the first subconstituent are different. For the pair $(U,W)$, if $X \in \Delta_1$ is a vertex adjacent to $U$ and $W$, then $X$ must contain $\langle e_1\rangle$, since $\langle e_2 \rangle \not\leq U$ and $\langle e_2 \rangle \not\leq W$. Hence, the number of common neighbours of $U$ and $W$ is equal to
\begin{align}
	\frac{q^n-1}{q-1}-4 = q^{n-1}+q^{n-2}+\ldots+q-3.\label{eq:cn1}
\end{align}
Now, let  $X \in \Delta_1$ be a common neighbour of $U^\prime$ and $W^\prime$ in the first subconstituent. Since $X\in \Delta_1$, we must have $V\cap X$ is non-trivial and $\alpha e_1+\beta e_2\in X$, for some $\alpha,\beta \in \mathbb{F}_q$ such that $(\alpha,\beta) \neq (0,0)$. Moreover, since $X\cap U^\prime$ and $X\cap W^\prime$ are both non-trivial, there exist non-zero elements $\gamma e_1+ \delta e_3,\gamma^\prime e_2+\delta^\prime e_3 \in X$, which are not necessarily distinct. 

First, assume that these elements are equal. Then, $\gamma e_1+ \delta e_3 =\gamma^\prime e_2+\delta^\prime e_3$ which implies $\gamma = \gamma^\prime = 0$, and $X\cap U^\prime = X\cap W^\prime = \langle e_3\rangle$. In this case, $X = \langle u,e_3 \rangle$, where $u\in V \setminus \{ \langle e_1\rangle,\langle e_2\rangle \}$ otherwise $X$ would not be in $\Delta_1$. Consequently, there are $q-1$ choices for such a subspace $X$.

Next, assume that $\gamma e_1+ \delta e_3 \neq\gamma^\prime e_2+\delta^\prime e_3$, or equivalently, $\langle e_3\rangle \not\leq X$. Hence, we have three pairwise linearly independent elements of $X$, which are $\alpha e_1+\beta e_2,\gamma e_1+ \delta e_3$, and $\gamma^\prime e_2+\delta^\prime e_3$. Since $W$ is a $2$-dimensional subspace, we have
\begin{align*}
	X = \langle \alpha e_1+\beta e_2 ,\gamma e_1+ \delta e_3\rangle = \langle \alpha e_1+\beta e_2, \gamma^\prime e_2+\delta^\prime e_3\rangle.
\end{align*}

We claim that $(\alpha,\beta)\neq (0,0)$.
If $\alpha = 0$, then $\langle e_2\rangle \in X$. So we have $X = \langle e_2,\gamma^\prime e_2+\delta^\prime e_3\rangle = W^\prime$ a contradiction, unless $\delta^\prime =0$ which also implies that $W$ is $1$-dimensional, also a contradiction. Hence, $\alpha \neq 0$. Similarly, one can use $\langle \alpha e_1+\beta e_2 ,\gamma e_1+ \delta e_3\rangle$ to show that $\beta \neq 0$.

We note that the subspace $\langle \alpha e_1+\beta e_2 ,\gamma e_1+ \delta e_3\rangle$ contains a unique $1$-dimensional subspace of $W^\prime = \langle e_2,e_3\rangle$, and similarly, $\langle \alpha e_1+\beta e_2, \gamma^\prime e_2+\delta^\prime e_3\rangle$ contains a unique $1$-dimensional subspace of $U^\prime = \langle e_1,e_3 \rangle$. Therefore, the number of choices for $X = \langle \alpha e_1+\beta e_2 ,\gamma e_1+ \delta e_3\rangle$ is exactly
\begin{align*}
	(q-1)^2.
\end{align*}
Overall, there are $q-1$ two-dimensional subspaces $X$ containing $\langle e_3\rangle$, and $(q-1)^2$ two-dimensional subspaces not containing $\langle e_3\rangle$. Therefore, the number of common neighbours of $U^\prime$ and $W^\prime$ is $q-1+(q-1)^2 = q(q-1).$

We conclude that the first subconstituent of $J_{q}(n,2)$ is not strongly regular, so again by Theorem~\ref{thm:triply-regular-equivalence} it cannot be triply regular and thus not triply transitive.
\begin{thm}
	For any $n\geq 4$ and a prime power $q$, the Grassmann graph $J_q(n,2)$ is not triply transitive.\label{thm:Grassmann}
\end{thm}

In the above theorem, we only need to show that the strongly-regular Grassmann graphs are not triply regular. It is worth noting that a consequence of the results in \cite{huang2025imperceptible} is that the Terwilliger algebra is strictly contained in the centralizer algebra.
\subsection{Paley graphs}
In has been shown in \cite{hanaki2023terwilliger} that for the Paley graph $P(q)$, where $q$ is a prime power, the centralizer algebra $\tilde{T}$ is larger than $T$, unless $q \in \{5,9\}$. The case where $q = 5$ is just the $5$-cycle, which is triangle free. The Paley graph $P(9)$ can be shown to be triply transitive since $\dim(\tilde{T}) = 15$.
\begin{thm}
	For any prime power $q\equiv 1\pmod 4$, the Paley graph $P(q)$ is triply transitive if and only if $q\in \{5,9\}$.\label{thm:Paley}
\end{thm}
\subsection{The Peisert graph $P^*(p^{2t})$}

Let $p\equiv 3\pmod 4$ be a prime number and $t\geq 1$ an integer. Let $q = p^{2t}$ and $\omega$ be a primitive element of $\mathbb{F}_{q}$. The \itbf{Peisert graph} $P^*(q)$ is the graph whose vertex set is $\mathbb{F}_{q}$, and $x,y\in \mathbb{F}_q$ are adjacent if and only if $y-x = \omega^i$, for $i\equiv 0\pmod 4$ or $i\equiv 1\pmod 4$. If $S = \langle \omega^4\rangle \cup \omega \langle \omega^4\rangle$, then it is not hard to see that $P^*(q)$ is in fact the Cayley graph $\operatorname{Cay}(\mathbb{F}_q,S)$.

Let us now recall some properties of the automorphism group of Peisert graphs. Define the maps $f_{a,b}: \mathbb{F}_q \to \mathbb{F}_q$, for $a\in \mathbb{F}_q^*$ and $b\in \mathbb{F}_q$, such that $v^{f_{a,b}} = av+b$, for any $v\in \mathbb{F}_q$. Let $\phi \in \Aut{\mathbb{F}_q/\mathbb{F}_p}$ be the Frobenius automorphism of $\mathbb{F}_q$. Note that the affine group of dimension $1$ is $\agl{1}{q} = \left\{ f_{a,b}: a\in \mathbb{F}_q^*,b\in \mathbb{F}_q \right\}$. It is well known that $ \Aut{\mathbb{F}_q/\mathbb{F}_p} = \langle \phi\rangle$ has order $2t$, and the affine semilinear group of dimension $1$ is $\agammal{1}{q} = \agl{1}{q}\rtimes \langle \phi\rangle$. For any $f_{a,b} \in \agl{1}{q}$ and $0\leq i\leq 2t-1$, we have
\begin{align*}
	v^{\phi^if_{a,b}} = \left(v^{p^i}\right)^{f_{a,b}} = av^{p^i}+b,\mbox{ for }v\in \mathbb{F}_q.
\end{align*}
Define the semilinear transformation $g_{a,b,i}: \mathbb{F}_q \to \mathbb{F}_q$ such that $v^{g_{a,b,i}} = \left(v^{p^i}\right)^{f_{a,b}} = av^{p^i}+b,$ for $v\in \mathbb{F}_q$. It is clear that $\agammal{1}{q} = \{ g_{a,b,i}: a\in \mathbb{F}_q^*,b\in \mathbb{F}_q,0\leq i \leq 2t-1 \}$.
Denote the subgroup of translation of $\agammal{1}{q}$ by 
\begin{align*}
	T = \left\{ g_{1,b,0} \in \agammal{1}{q}: b\in \mathbb{F}_q \right\}.
\end{align*}
Denote the scalar multiplication by $\omega$ with $\varphi$. Define the group
\begin{align}
	G_q = T \langle \varphi^4,\phi\varphi\rangle.\label{eq:aut}
\end{align}

\begin{thm}\cite{peisert2001all}
	If $q\neq 3^2,7^2,3^4,23^2$, then $\Aut{P^*(q)} = G_q$. In particular, $|\Aut{P^*(p^{2t})}| = \tfrac{p^{2t}(p^{2t}-1)2t}{4}$. Moreover, the automorphism groups of $P^*(3^2),P^*(7^2),P^*(3^4), \mbox{ and }P^*(23^2)$ are respectively $\sym{3}\wr \sym{2}$, $7^2:(3\times \sln{2}{3})$, $3^4:(\sln{2}{5}:2^2)$, and $23^2:(11\times(3:Q_8))$.\label{thm:Aut-Peisert}
\end{thm}

The Peisert graph $P^*(9)$ is isomorphic to the Paley graph $P(9)$. We have seen in the previous section that $P(9)$ is triply transitive. Using \verb*|Sagemath|, the dimensions of the centralizer algebras of $P^*(49)$, $P^*(81)$, and $P^*(529)$ are respectively $45$, $31$, and $1061$. We conclude that these graphs cannot be triply transitive.

	Now, we assume that $q = p^{2t}\neq 3^2,7^2,3^4,23^2$ and consider the Peisert graph $P^*(q)$. By Theorem~\ref{thm:Aut-Peisert}, we know that $\Aut{P^*(q)} = G_q$.
We note that $\Delta_0 = \{0\}$, $\Delta_1 = \langle \omega^4\rangle \cup \omega\langle \omega^4\rangle$, and $\Delta_2 = \{ \omega^i : i\not \equiv 0,1\pmod 4 \}$. The stabilizer of $0$ in $G_q$ is the group $H = \langle \varphi^4,\phi\varphi \rangle$. Moreover, for any $0\leq i \leq q-1$, we have
\begin{align*}
	(\omega^i)^{\phi\varphi\phi^{-1}} =\omega^{(ip+1)p^{2t-1}} = \omega^{i+p^{2t-1}} = \left(\omega^{i}\right)^{\varphi^{p^{2t-1}}}.
\end{align*}
Consequently, $\phi\varphi\phi^{-1} = \varphi^{p^{2t-1}}$, and so 
\begin{align*}
	\left(\phi\varphi\right) \varphi^4 \left(\phi\varphi\right)^{-1} = \phi \varphi^4 \phi^{-1} = \varphi^{4p^{2t-1}}.
\end{align*}
We deduce that $\langle \varphi^4\rangle \trianglelefteq H$, and $H = \langle \varphi^4\rangle \langle \phi \varphi\rangle$. As $G$ is an affine group, $G = T\rtimes H$, where $H$ is the full point stabilizer. Therefore, 
\begin{align*}
	\tfrac{(p^{2t}-1)t}{2} = |H| = \tfrac{|\langle \varphi^4 \rangle|\ |\langle \phi\varphi \rangle|}{|\langle \varphi^4 \rangle \cap \langle \phi \varphi \rangle|}.
\end{align*}
We conclude that $\tfrac{ |\langle \phi\varphi \rangle|}{|\langle \varphi^4 \rangle \cap \langle \phi \varphi \rangle|} = 2t.$ In particular, $2t\mid o(\phi\varphi)$.

\begin{prop}
	We have $\langle \varphi^4\rangle \cap \langle \phi \varphi\rangle$ is trivial.\label{prop:trivial}
\end{prop}
\begin{proof}
	If $\psi \in\langle \varphi^4\rangle \cap \langle \phi \varphi\rangle$, then there exists $0\leq r \leq \tfrac{p^{2t}-1}{4}-1$ and a non-negative integer $s$ such that $\psi = \varphi^{4r}=(\phi\varphi)^s$. So, for $0\leq i\leq q-2$, we have $\left(\omega^i\right)^{\varphi^{4r}}=\left(\omega^i\right)^{(\phi\varphi)^s}$. Hence,
	\begin{align*}
		\omega^{i+4r}=\omega^{ip^s+p^{s-1}+\ldots+p+1}
	\end{align*}
	for all $0\leq i\leq q-2$. Equivalently, we have
	\begin{align}
		\left(1-p^s\right)i+4r\equiv p^{s-1}+\ldots+p+1 \pmod{p^{2t}-1},\label{eq:intersection}
	\end{align}
	 for all $0\leq i\leq q-2$. If $i = p-1$, then \eqref{eq:intersection} becomes
	 \begin{align*}
	 	4r &\equiv p^{s+1} - p^s -p +1 +p^{s-1}+ \ldots + p + 1 \pmod{p^{2t}-1}\\
	 	&\equiv p^{s+1} - p^s + p^{s-1}+ \ldots+p^2+2 \pmod{p^{2t}-1}.
	 \end{align*}
	 Since this equation has a solution in $r$ and $s$, we must have
	 \begin{align}
	 	\gcd\left(4,p^{2t}-1\right) \Big| \left(p^{s+1}-p^s+p^{s-1}+\ldots+p^2+2\right).\label{eq:division}
	 \end{align}
	 
	We have $\gcd\left(4,p^{2t}-1\right) = 4$, and since $p\equiv -1 \pmod 4$, we also have
	 \begin{align*}
	 	p^{s+1} - p^s + p^{s-1}+ \ldots+p^2+2 \equiv (-1)^{s+1}-(-1)^s+(-1)^{s-1}+\ldots+(-1)^2 +2  \pmod{4}. 
	 \end{align*}
	 If $s$ is even, then $(-1)^{s+1}-(-1)^s+(-1)^{s-1}+\ldots+(-1)^2 +2 \equiv 0 \pmod{4}$, and if $s$ is odd, then $(-1)^{s+1}-(-1)^s+(-1)^{s-1}+\ldots+(-1)^2 +2 \equiv 1 \pmod{4}$. In this case, we conclude that $s$ must be even. In particular, we have 
	 \begin{align}
	 	r\equiv \frac{p^{s+1} - p^s + p^{s-1}+ \ldots+p^2+2}{4} \pmod{\tfrac{p^{2t}-1}{4}}\label{eq:Peisert-minus}
	 \end{align}
	 and there are exactly $4$ solutions modulo $p^{2t}-1$.	 
	 
	If $i = 0$, then \eqref{eq:intersection} becomes
	 \begin{align*}
	 	4r\equiv p^{s-1}+\ldots+p+1 \pmod{p^{2t}-1}.
	 \end{align*}
	 The above equation has a solution in $r\in\mathbb{F}_q$ if and only if $s$ is even. Hence, we have 
	 \begin{align}
	 	r\equiv \frac{p^s-1}{4(p-1)} \pmod{\tfrac{p^{2t}-1}{4}}.\label{eq:Peisert-1}
	 \end{align}
	From \eqref{eq:Peisert-minus} and \eqref{eq:Peisert-1}, we conclude that 
	\begin{align*}
		\frac{p^{s-1}+\ldots+p+1}{4}\equiv \frac{p^{s+1} - p^s + p^{s-1}+ \ldots+p^2+2}{4} \pmod{\tfrac{p^{2t}-1}{4}}
	\end{align*}
	or equivalently
	\begin{align*}
		\frac{p^{s+1} - p^s -p+1}{4} \equiv 0 \pmod{\tfrac{p^{2t}-1}{4}}.
	\end{align*}
	We deduce that $p^{2t}-1 \mid (p-1)(p^s-1)$ which can only happen if $s = 0$. Therefore, $r = 0$, and $\psi$ is trivial. This completes the proof.
\end{proof}
We conclude from Proposition~\ref{prop:trivial} that the order of $\phi\varphi$ is exactly $2t$.
Now, we determine the number of orbitals of $H$ in $\Delta_1$. Note that $H$ acts transitively on $\Delta_1$, $|\Delta_1| = \tfrac{p^{2t}-1}{2}$, and a stabilizer of a point of $\Delta_1$ has order 
\begin{align*}
	\frac{|H|}{|\Delta_1|} = \tfrac{\tfrac{(p^{2t}-1)t}{2}}{\tfrac{(p^{2t}-1)}{2}}=t. 
\end{align*}

\begin{lem}
	We have
	\begin{align*}
		H_1 &= \left\{ \varphi^{-\tfrac{p^{2k}-1}{p^{2k}(p-1)}} (\phi\varphi)^{2k} : 0\leq k \leq 2t-1 \mbox{ is even} \right\}.
	\end{align*}
\end{lem}
\begin{proof}
	If $\varphi^{4r} (\phi\varphi)^s\in H$ fixes $1$, then we have
	\begin{align*}
		1 &=1^{\varphi^{4r} (\phi\varphi)^s}\\
		&=\left(\omega^{4r}\right)^{(\phi\varphi)^s}\\
		&=\left(\omega^{4rp+1}\right)^{(\phi\varphi)^{s-1}}\\
		&\ \ \vdots\\
		&=\omega^{4rp^s+p^{s-1}+\ldots+p+1}.
	\end{align*}
	The above equation holds if and only if $4rp^s \equiv -\left(p^{s-1}+\ldots+p+1\right) \pmod{p^{2t}-1}$. The latter has a solution in $r$ if and only if $\gcd(4p^s,p^{2t}-1) = 4$ divides $ \left(p^{s-1}+\ldots+p+1\right)$. Since $p \equiv -1 \pmod{4}$,  $4\mid \left(p^{s-1}+\ldots+p+1\right)$ if and only if $s$ is even. 
	For $s$ even, we have
	\begin{align*}
		r \equiv -\frac{p^s-1}{4p^s(p-1)} \pmod{\tfrac{p^{2t}-1}{4}}.
	\end{align*}
	Hence, if $s = 2k$, then we have
	\begin{align*}
		\varphi^{4r} (\phi\varphi)^s = \varphi^{-\tfrac{p^{2k}-1}{p^{2k}(p-1)}} (\phi\varphi)^{2k}.
	\end{align*}
	This completes the proof.
\end{proof}
Let us know determine the multiplication in $H_1$ and show that the latter must be cyclic. For any $\ell$ and $k$ such that $0\leq 2\ell,2k \leq 2t-1$, we have
\begin{align*}
	\left(\varphi^{-\tfrac{p^{2k}-1}{p^{2k}(p-1)}} (\phi\varphi)^{2k}\right)\left(\varphi^{-\tfrac{p^{2\ell}-1}{p^{2\ell}(p-1)}} (\phi\varphi)^{2\ell}\right)
	&= \varphi^{-\left(\tfrac{p^{2k}-1}{(p-1)p^{2k}}+p^{2k(2t-1)}\tfrac{p^{2\ell}-1}{(p-1)p^{2\ell}}\right)} \left(\phi\varphi\right)^{2k+2\ell}\\
	&= \varphi^{-\left(\tfrac{p^{2k+2\ell}-p^{2\ell}+p^{(2t)(2k)}(p^{2\ell}-1)}{(p-1)p^{2k+2\ell}}\right)} \left(\phi\varphi\right)^{2k+2\ell}\\
	&= \varphi^{-\left(\tfrac{p^{2k+2\ell}-p^{2\ell}+(p^{2\ell}-1)}{(p-1)p^{2k+2\ell}}\right)} \left(\phi\varphi\right)^{2k+2\ell}\\
	&= \varphi^{-\left(\tfrac{p^{2k+2\ell}-1}{(p-1)p^{2k+2\ell}}\right)} \left(\phi\varphi\right)^{2k+2\ell}.
\end{align*}
Therefore, we have
\begin{align*}
	\left(\varphi^{-\tfrac{p^{2k}-1}{p^{2k}(p-1)}} (\phi\varphi)^{2k}\right)\left(\varphi^{-\tfrac{p^{2\ell}-1}{p^{2\ell}(p-1)}} (\phi\varphi)^{2\ell}\right)
	&=
	\varphi^{-\left(\tfrac{p^{2k+2\ell}-1}{(p-1)p^{2k+2\ell}}\right)} \left(\phi\varphi\right)^{2k+2\ell}.
\end{align*}
From this, we deduce that $H_1$ is cyclic, and so
\begin{align*}
	H_1 = \left\langle \varphi^{-\tfrac{p^2-1}{p^2(p-1)}}\left(\phi\varphi\right)^2 \right\rangle.
\end{align*}
Next, we determine the number of fixed points of the permutation $\varphi^{-\tfrac{p^2-1}{p^2(p-1)}}\left(\phi\varphi\right)^2$. For any $0\leq i\leq p^{2t}-2$, if $\omega^i$ is fixed by $\varphi^{-\tfrac{p^2-1}{p^2(p-1)}}\left(\phi\varphi\right)^2$, then we have
\begin{align*}
	\omega^i &= \left(\omega^i\right)^{\varphi^{-\tfrac{p^2-1}{p^2(p-1)}}\left(\phi\varphi\right)^2}\\
	&= \omega^{\left(i-\tfrac{p^2-1}{p^2(p-1)}\right)p^2+p+1}\\
	&= \omega^{\left(ip^2-(p+1)\right)+p+1}\\
	&= \omega^{ip^2}.
\end{align*}
Therefore, $\omega^i$ must lie in $\mathbb{F}_{p^2}$. Hence, $\varphi^{-\tfrac{p^2-1}{p^2(p-1)}}\left(\phi\varphi\right)^2$ fixes every element of $\mathbb{F}_{p^2}^*$. In particular, $H_1$ has at least $p^2$ orbits.

Now, we are ready to prove the main result.
\begin{lem}
	We have $\dim(\tilde{T})>15$.
\end{lem}
\begin{proof}
	Consider the Peisert graph $P^*(p^{2t})$, where $p\equiv 3 \pmod 4$. Recall the full automorphism group of $P^*(p^{2t})$ is the group $G_q$ given in \eqref{eq:aut}, the stabilizer of $0 \in \mathbb{F}_{p^{2t}}$ is the group $H = \langle \varphi^4,\phi\varphi\rangle$, and $H_1 =\langle \varphi^{-\tfrac{p^2-1}{p^2(p-1)}}\left(\phi\varphi\right)^2 \rangle$.
	
	 If $t=1$, then $|H_1| = t=  1$. As $P^*(p^{2t})$ is self-complementary, the group $H_\omega$, for any $\omega \in \Delta_2$, has the same number of orbits as $H_1$. Hence, we have 
	 \begin{align*}
	 	\dim(\tilde{T})\geq 5+ 2(p^2-1)\geq 60
	 \end{align*} 
	 whenever $p\geq 5$ (the case $p=3$ is an exception since the Peisert graph with this parameter is the same as a Paley graph on $9$ vertices).
	
	Now, assume that $t\geq 2$. Since the group $H_1$ fixes every element of $\mathbb{F}_{p^2}^*$, we know that $H_1$ has at least $p^2-1$ orbits on the vertices of the Peisert graph $P^*(p^{2t})$. As the Peisert graph is a self-complementary graph, the group $H_x$ for $x\in \Delta_2$ also has at least $p^2-1$ orbits. Therefore, we have
	\begin{align*}
		\dim(\tilde{T})\geq 5+(p^2-1)+(p^2-1)>15.
	\end{align*} 
	This completes the proof.
\end{proof}
\begin{thm}
	The Peisert graph $P^*(p^{2t})$ is triply transitive only when $(p,t) = (3,2)$.\label{thm:Peisert}
\end{thm}

\subsection{The $n\times n$ grid}\label{subsect:grid}
The $n\times n$ grid is isomorphic to the Hamming graph $H(2,n)$. It is clear that this graph is the Latin square graph $\operatorname{LS}_{2}(n)$, and parameters are
\begin{align*}
	\left(n^2,2(n-1),n-2,2\right).
\end{align*}

In \cite{levstein2006terwilliger}, it was shown that the dimension of the Terwilliger algebra is $15$ if $n\geq 3$ and $10$ if $n = 2$. We verify using \verb*|Sagemath| that $H(2,2)$ is triply transitive. For $n\geq 3$, the graph $H(2,n)$ is triply regular since $\dim(T_0) = 15$. In what remains, we will show that $\dim(\tilde{T}) = 15$. 

Recall that the automorphism group of $H(2,n) = K_n\square K_n$ is $\sym{n} \wr \sym{2}$ and the stabilizer of a vertex is isomorphic to $\sym{n-1}\wr \sym{2}$. 

\begin{thm}
	For any $n\geq 3$, the $n\times n$ grid is triply transitive.\label{thm:grid}
\end{thm}
\begin{proof}
	Since $\Gamma$ and its complement do not admit triangles, we know that $\dim(T_0) = 15$. It remains to show that the centralizer algebra $\tilde{T}$ has dimension $15$. The automorphism group of $H(2,n)$ is $\sym{n}\wr \sym{2}$ acting on $\{1,2,\ldots,n\}\times \{1,2,\ldots,n\}$. Let $\varphi \in \Aut{\sym{n}\times \sym{n}}$ be the involution with the property $(i,j)^\varphi = (j,i)$, for $(i,j) \in \{1,2,\ldots,n\}\times \{1,2,\ldots,n\}$. The elements of $\sym{n}\wr \sym{2}$ are of the form $(\sigma,\tau;\varphi^k)$ such that $\sigma,\tau \in \sym{n}$, and $0\leq k\leq 1$. For $(i,j) \in \{1,2,\ldots,n\} \times \{1,2,\ldots,n\}$ and $(\sigma,\tau;\varphi^k) \in \sym{n}\wr \sym{2}$, we have 
	\begin{align*}
		(i,j)^{(\sigma,\tau;\varphi^k)} =
		\begin{cases}
			(i^\sigma,j^\tau) &\mbox{ if }k = 0\\
			(j^\tau,i^\sigma) & \mbox{ if } k = 1.
		\end{cases}
	\end{align*}
	The stabilizer of $(n,n)$ in $\sym{n}\wr \sym{2}$ is $\sym{n-1}\wr \sym{2}$. Let us count the number of orbitals of $\sym{n-1}\wr \sym{2}$. Note that 
	\begin{align*}
		\Delta_0 &= \{(n,n)\}, \\
		\Delta_1 &= \{ (i,j): i=n \mbox{ and } 1\leq j\leq n-1, \mbox{ or } j=n \mbox{ and } 1\leq i\leq n-1 \}\\
		\Delta_2 &= \{ (i,j): i,j \neq n \}.
	\end{align*}
	 As we have seen before, this is equivalent to counting certain orbitals.
	\begin{itemize}
		\item As $H(2,n)$ is a rank $3$ graph, we know that $\sym{n-1}\wr \sym{2}$ has $3$ orbits on the vertices.
		\item The vertex $(n,n-1)$ is adjacent to $(n,n)$ in $H(2,n)$. The stabilizer of the vertex $(n,n-1)$ in $\sym{n-1}\wr \sym{2}$ is the subgroup $\sym{n-1}\times \sym{n-2}$. There is a unique orbit of this subgroup in $\Delta_0$. 
		
		The set $\Delta_1$ splits into three orbits. The first orbit is $\{(n,n-1)\}$. The other two orbits are the sets
		\begin{align*}
			\left\{ (n,j): 1\leq j\leq n-2 \right\} \mbox{ and } \left\{ (i,n): 1\leq i\leq n-1 \right\}.
		\end{align*}
		
		The set $\Delta_2$ splits into two orbits of $\sym{n-1} \times \sym{n-2}$. These orbits are
		\begin{align*}
			\left\{ (i,j): 1\leq i,j\leq n-2 \right\} \cup \left\{ (n-1,j): 1\leq j\leq n-2 \right\}
		\end{align*}
		and 
		\begin{align*}
			\left\{ (i,n-1): 1\leq i\leq n-2 \right\}.
		\end{align*}
		\item Now, consider the vertex $(n-1,n-1) \in \Delta_2$. The stabilizer of this vertex in $\sym{n-1}\wr \sym{2}$ is the subgroup $\sym{n-2}\wr \sym{2}$. $\Delta_0$ is an orbit of this subgroup. The set $\Delta_1$ splits into two orbits of $\sym{n-2}\wr \sym{2}$, which are the sets 
		\begin{align*}
			\left\{ (i,n):1\leq i\leq n-2 \right\}\cup \left\{ (n,j): 1\leq j\leq n-2 \right\} \mbox{ and } \{ (n-1,n),(n,n-1) \}.
		\end{align*}
		The set $\Delta_2$ clearly splits into $3$ orbits of $\sym{n-2}\wr \sym{2}$ since the subgraph of $H(2,n)$ induced by $\Delta_2$ is the rank $3$ graph $H(2,n-1)$. Hence, the stabilizer $\sym{n-2}\wr \sym{2}$ of $(n-1,n-1)$ in $H(2,n-1)$ has $3$ orbits.
	\end{itemize} 
	We conclude that $\dim(\tilde{T}) = 3+6+6=15$, so $H(2,n)$ is triply transitive.
\end{proof}

\subsection{The graph $E_{6,1}(q)$}

The first subconstituent of the graph $E_{6,1}(q)$ (see  \cite[\S4.9]{brouwer2022strongly}) is the $q$-clique extension of the collinearity graph of $\vd_{5,5}(q)$ defined in \cite[\S3.3.3]{brouwer2022strongly}, which is a primitive rank $3$ graph. By Corollary~\ref{cor:m-extension}, we deduce that this first subconstituent of $E_{6,1}(q)$ is not strongly regular, and therefore $E_{6,1}(q)$ cannot be triply transitive.
\begin{thm}
	The graph $E_{6,1}(q)$ is not triply transitive for any prime power $q$.\label{thm:e61}
\end{thm}

\subsection{Collinearity graphs}

\subsubsection{Classical polar spaces}
Let $\mathcal{P}$ be a finite classical polar space with rank $n\geq 2$ on a vector space $V$. Let $\Gamma$ be the collinearity graph of $\mathcal{P}$. By \cite[pg~37]{brouwer2022strongly}, the collinearity graph of $\mathcal{P}$ is locally the $q$-clique extension of the collinearity graph of a polar space of rank $n-1$. By Corollary~\ref{cor:m-extension}, no first subconstituent of $\Gamma$ is strongly regular for $n\geq 3$. Hence, $\Gamma$ is not triply transitive whenever $n\geq 3$.

\begin{thm}
	If $\mathcal{P}$ is a classical polar space of rank $n\geq 3$, then its collinearity graph is not triply transitive.\label{thm:polar-graph-rank-3}
\end{thm}

Now, assume that $n = 2$. We distinguish the cases depending on the form associated with the polar space $\mathcal{P}$. 
We summarize the triple transitivity status of the collinearity graph of rank $2$ polar spaces in the next table.
\begin{table}[H]
	\centering
	\begin{tabular}{ccccc}
		\hline
		Row&$\mathcal{P}$&Order& Triply transitive & Comments\\
		\hline \hline
		1&$\operatorname{O}_4^+(q)$&$(q,1)$&Yes&$(q+1)\times (q+1)$ grid\\
		2&$\operatorname{Sp}_4(q)$&$(q,q)$&No& $\mathrm{q}_{11}^1,\mathrm{q}_{22}^2>0$\\
		3&$\operatorname{O}_5(q)$&$(q,q)$&No&$\mathrm{q}_{11}^1,\mathrm{q}_{22}^2>0$\\
		4&$\operatorname{O}_6^-(q)$&$(q,q^2)$&?&Smith graph(thus triply regular)\\
		5&$\operatorname{U}_4(\sqrt{q})$&$(q,q)$&No&$\mathrm{q}_{11}^1,\mathrm{q}_{22}^2>0$\\
		6&$\operatorname{U}_5(\sqrt{q})$&$(q,q^3)$&No&$\mathrm{q}_{11}^1,\mathrm{q}_{22}^2>0$\\
	\end{tabular}
	\caption{Status of triple transitivity for the collinearity graphs of classical polar spaces of rank $2$.}\label{tab:rank2}
\end{table}

In the next theorem, we identify the collinearity graph of classical polar spaces of rank $2$ that are not triply transitive.
\begin{thm}
	The collinearity graph of the polar spaces in Table~\ref{tab2} are:
	\begin{enumerate}[(a)]
		\item triply transitive for the polar space in Row~1,
		\item not triply transitive for the polar spaces in Rows~2,3,5,6.
	\end{enumerate}
	\label{thm:polar-graph-rank-2}
\end{thm}
\begin{proof}
	Let $\Gamma$ be the collinearity graph of the polar in Table~\ref{tab:rank2}. 
	
	\begin{itemize}
		\item For Row~1, the collinearity graph $\Gamma$ is isomorphic to the $(q+1)\times (q+1)$ grid, which is triply transitive by Theorem~\ref{thm:grid}. 
		\item For Row~2, the number of vertices of $\Gamma$ is $\tfrac{q^4-1}{q-1}$. By a result of Ljunggren \cite{ljunggren1943noen}, the number $\tfrac{q^4-1}{q-1}$ is a perfect square only when $q=7$. In this case, the collinearity graph $\Gamma$ has parameters $(400, 56, 6, 8)$, which is not a Latin square or negative Latin square parameters. Hence, the collinearity graphs $\Gamma$ do not have Latin square or negative Latin square parameters. By Lemma~\ref{lem:srg-not-triply-regular}, and the comments in Table~\ref{tab2}, we deduce $\Gamma$ is not triply transitive.
		\item For Row~3, the graph $\Gamma$ has the same parameters as the collinearity graph of $\operatorname{Sp}_4(q)$, so it cannot have Latin square of negative Latin square parameters. Similarly, Lemma~\ref{lem:srg-not-triply-regular} implies that $\Gamma$ is not triply transitive.
		
		\item Next, we consider Row~5. The graph $\Gamma$ has parameters
		\begin{align}
			\left((q^2+1)(q^3+1),q^2(q+1),q^2-1,q+1\right).\label{eq:U4}
		\end{align}
		Note that \eqref{eq:U4} are not the parameters of a grid and also not of the form $(r^2(r+3)^2,r^3+3r^2+r,0,r^2+r)$, for some integer $r\geq 2$. 
		
		If \eqref{eq:U4} are the parameters of a Latin square or negative Latin square graph, then by Lemma~\ref{lem:rank3-LS-nLS} and the fact that $q^2-1\neq q+1$, we deduce that at least one of the subconstituents of $\Gamma$ is not strongly regular. Hence, $\Gamma$ is not triply transitive. 
		
		If \eqref{eq:U4} are not the parameters of a Latin square or negative Latin square graph, then using the fact that $\mathrm{q}_{11}^1,\mathrm{q}_{22}^2>0$, we deduce by Lemma~\ref{lem:srg-not-triply-regular} that $\Gamma$ is not triply transitive.
		
		\item For Row~6, the graph $\Gamma$ has parameters
		\begin{align*}
			\left((q^2+1)(q^5+1),q^2(q^3+1),q^2-1,q^3+1\right).
		\end{align*}
		Using the same argument as Row~5, we conclude that $\Gamma$ is not triply transitive.
	\end{itemize}
	This completes the proof.
\end{proof}
Consider the collinearity graph $\Gamma$ of the polar space in Row~4 of Table~\ref{tab:rank2}. The number of vertices of $\Gamma$ is $(q+1)(q^3+1)$. Note that $(q+1)(q^3+1) = (q+1)^2(q^2-q+1)$ and $(q-1)^2<q^2-q+1< q^2 $. Hence, the number of vertices of $\Gamma$ is not a square. Consequently, $\Gamma$ does not have a Latin square or negative Latin square parameters. However, $\Gamma$ is a Smith graph, and hence, $\mathrm{q}_{ii}^i = 0$ for some $i\in \{1,2\}$, so we cannot use Lemma~\ref{lem:srg-not-triply-regular}. Therefore, the status of this graph is open.
Based on computations for small values of $q \in \{2,3,4,5\}$, we conjecture the following for the polar space in Row~4 of Table~\ref{tab2}.
\begin{conj}
	The collinearity graph of polar space $O_6^-(q)$ is triply transitive.
\end{conj}

\subsubsection{Dual of a finite classical polar space of rank $2$}

Let $\Gamma$ be the collinearity graph of the dual a finite classical polar space $\mathcal{P}$ of rank $2$. Then, $\mathcal{P}$ corresponds to a generalized quadrangle. By \cite[Theorem~2.2.19]{brouwer2022strongly}, $\Gamma$ is the collinearity graph of a generalized quadrangle. These graphs arise from the polar spaces in Table~\ref{tab:rank2-quadrangle}.

\begin{table}[b]
	\centering
	\begin{tabular}{ccccc}
		\hline
		Row&$\mathcal{P}$&Order& Triply transitive & Comments\\
		\hline \hline
		1&$\operatorname{O}_4^+(q)$&$(q,1)$&Yes&$K_{(q+1),(q+1)}$ (imprimitive)\\
		2&$\operatorname{Sp}_4(q)$&$(q,q)$&No& $\mathrm{q}_{11}^1,\mathrm{q}_{22}^2>0$\\
		3&$\operatorname{O}_5(q)$&$(q,q)$&No&$\mathrm{q}_{11}^1,\mathrm{q}_{22}^2>0$\\
		4&$\operatorname{O}_6^-(q)$&$(q,q^2)$&No&$\mathrm{q}_{11}^1,\mathrm{q}_{22}^2>0$\\
		5&$\operatorname{U}_4(\sqrt{q})$&$(q,q)$&No&$\mathrm{q}_{11}^1,\mathrm{q}_{22}^2>0$\\
		6&$\operatorname{U}_5(\sqrt{q})$&$(q,q^3)$&No&$\mathrm{q}_{11}^1,\mathrm{q}_{22}^2>0$\\
	\end{tabular}
	\caption{Status of triple transitivity for the collinearity graphs of the dual of the classical polar spaces of rank $2$.}\label{tab:rank2-quadrangle}
\end{table}

Similar to the collinearity graphs of polar spaces of rank $2$, the number of vertices of the collinearity graphs of generalized quadrangles arising from Rows~2-6 in Table~\ref{tab:rank2-quadrangle} are not square. Consequently, the graph $\Gamma$ for these rows do not have a Latin square or negative Latin square parameters. By the comments in Table~\ref{tab:rank2-quadrangle} and Lemma~\ref{lem:srg-not-triply-regular}, we deduce the following corollary.
\begin{thm}
	The collinearity graphs of the generalized quadrangle in Rows~2-6 of Table~\ref{tab:rank2-quadrangle} are not triply transitive.\label{thm:polar-dual-graph-rank-2}
\end{thm}

\subsection{Affine polar graph}\label{subs:affine-polar-graph}
Consider the affine polar graph $VO^\varepsilon_{2m}(q)$ defined in \cite[\S3.3]{brouwer2022strongly}, where $m\geq 1$ is an integer, $q$ is a prime power, and $\varepsilon \in \{\pm 1\}$. 
We note that the graph $\vo_{2m}^\varepsilon(q)$ has parameters equal to 
\begin{align*}
	\left( q^{2m},(q^m-\varepsilon)(q^{m-1}+\varepsilon),q(q^{m-1}-\varepsilon)(q^{m-2}+\varepsilon)+q-2,q^{m-1}(q^{m-1}+\varepsilon) \right).
\end{align*}

\begin{rmk}
	The graph $\vo_{2m}^\varepsilon(q)$ has a Latin square parameter and negative Latin square parameters if and only if $(q,\varepsilon)=(2,1)$ and $(q,\varepsilon) = (2,-1)$, respectively.\label{rmk:affine-polar-graph}
\end{rmk}

Let us first consider the cases where $m = 1$. The graph $VO_2^-(q)$ is a coclique, so the automorphism group has rank $2$. The graph $VO^+_{2}(q)$ is the $q\times q$ grid, so it is triply transitive.

Now assume that $m\geq 2$. As explained at the beginning of \S3.3.3 in \cite{brouwer2022strongly}, the first subconstituent of $VO_{2m}^\varepsilon(q)$ is the $(q-1)$-clique extension of the graph $\Gamma(O^\varepsilon_{2m}(q))$, where $\Gamma(O^\varepsilon_{2m}(q))$ is the collinearity graph of the orthogonal polar space defined in \cite[\S2.6.3]{brouwer2022strongly}. Since the latter is primitive, by Corollary~\ref{cor:m-extension}, we deduce that $VO_{2m}^\varepsilon(q)$ is not triply transitive for $q>2$. 

\begin{thm}
	The graph $\vo_{2m}^\varepsilon(q)$ is not triply transitive for $\varepsilon = \pm $, $q>2$, and $m\geq 2$.\label{thm:affine-polar-graph}
\end{thm}

The graph $VO^\varepsilon_{2m}(2)$ in triply regular by \cite[Pg~25]{brouwer2022strongly} and \cite[Proposition~3.6.1]{brouwer2022strongly}. 
In \cite{brouwer1990graphs} and \cite[Proposition~3.6.1]{brouwer2022strongly}, Brouwer and Shult showed that the first and second subconstituents of $VO_{2m}^\varepsilon(2)$ are respectively the graphs $\Gamma(O^\varepsilon_{2m}(2))$ and $NO_{2m}^\varepsilon(2)$. Since these two graphs are also rank $3$ graphs, we deduce that the block dimension decomposition of $\tilde{T}$ for $VO_{2m}^\varepsilon(2)$ is of the form
\begin{align*}
	\begin{bmatrix}
		1 & 1 & 1\\
		1 & 3 & t\\
		1 & t & 3
	\end{bmatrix}
\end{align*}
for some integer $t\geq 1$. As the subconstituents of $VO_{2m}^\varepsilon(2)$ are both strongly regular, we deduce from Theorem~\ref{thm:triply-regular-equivalence} that $T_0 = T$. Consequently, we must have $t\geq 2$. We conjecture that in fact $t = 2$.

\begin{conj}
	The affine polar graph $VO^\varepsilon_{2m}(2)$ is triply transitive.
\end{conj}

\subsection{Affine half spin graph}

The affine half spin graph $\vd_{5,5}(q)$ is defined in \cite[\S3.3.3]{brouwer2022strongly}. The first subconstituents of the graph $\vd_{5,5}(q)$ are isomorphic to the $(q-1)$-clique extension of the halved graph of the collinearity graph of the orthogonal polar space $O^{+}_{10}(q)$. Hence, for $q>2$, the graph $\vd_{5,5}(q)$ is not triply transitive. 

It remains to check the graph $\vd_{5,5}(2)$. This graph has parameters
$
(65536,2295,310,72).
$
Hence, the graph $\vd_{5,5}(2)$ has Latin square parameters. Using Lemma~\ref{lem:rank3-LS-nLS}, we deduce that $\vd_{5,5}(2)$ is not triply transitive.
\begin{thm}
	The affine half spin graph $\vd_{5,5}(q)$ is not triply transitive for any prime power $q$.\label{thm:affine-half-spin}
\end{thm}

\subsection{The graph $\vsz(q)$, for $q=2^{2e-1}$}
Consider the  graph $\vsz(q)$, for $q=2^{2e-1}$ defined in \cite[\S3.3]{brouwer2022strongly}. The graph $\vsz(q)$ is a rank $3$ graph with parameters equal to those of $\vo_4^-(q)$. Hence, the parameters of $\vsz(q)$ are
\begin{align*}
	\left(q^4,(q^2+1)(q-1),q-2,q(q-1)\right).
\end{align*}
Hence, $\vsz(q)$ has negative Latin square parameters. The parameters of its complement are
\begin{align*}
	\left(q^4,q^4 - q^3 + q^2 - q,q^{4} - 2 \, q^{3} + 3 \, q^{2} - 3 \, q,q^{4} - 2 \, q^{3} + 2 \, q^{2} - q\right). 
\end{align*}

If $\vsz(q)$ is triply transitive, then by Lemma~\ref{lem:rank3-LS-nLS}, it has parameters
\begin{align*}
	(4t^2,t(2t+1),t(t+1),t(t+1)).
\end{align*}
Therefore, we must have $q-2=q(q-1)$, which implies $(q-1)^2+1=0$, a contradiction, or $$q^{4} - 2 \, q^{3} + 3 \, q^{2} - 3 \, q=q^{4} - 2 \, q^{3} + 2 \, q^{2} - q$$
implying that $q = 2$. The graph $\vsz(2)$ in this case is triply transitive and is isomorphic to $\vo_4^-(2)$.
\begin{thm}
	The graph $\vsz(q)$ is triply transitive for $q = 2^{2e-1}$ if and only if $e = 1$.\label{thm:ST}
\end{thm}

\subsection{Alternating forms graphs}
The alternating forms graph $\alt{5,q}$ is the graph whose vertex set consists of the $5\times 5$ skew-symmetric matrices over $\mathbb{F}_q$, where two matrices are adjacent if and only if the rank of their difference is $2$. The graph $\alt{5,q}$ is a rank $3$ graph with parameters
\begin{align*}
	\left(q^{10},(q^2+1)(q^5-1),q^5+q^4-q^2-2,q^2(q^2+1)\right).
\end{align*}
Therefore, $\alt{5,q}$ has Latin square parameters. The parameters of its complement is
\begin{align*}
	\left(q^{10}, q^{10} - q^{7} - q^{5} + q^{2}, q^{10} - 2 \, q^{7} - 2 \, q^{5} + q^{4} + 3 \, q^{2}, q^{10} - 2 \, q^{7} - q^{5} + q^{4} + q^{2}\right).
\end{align*}

If $\alt{5,q}$ is triply transitive, then its parameters must be of the form
\begin{align*}
	(4t^2,t(2t-1),t(t-1),t(t-1))
\end{align*}
by Lemma~\ref{lem:rank3-LS-nLS}, for some $t\geq 2$. Consequently, we have $2t = q^5$ implying that $q$ is a power of $2$. Moreover, we have $t = q^2+1$. Therefore, we have $q^5 = 2t = 2(q^2+1)$, which has no solution in $q$ as the latter must be a power of $2$.

For the complement, we must have
\begin{align*}
	q^{10} - 2 \, q^{7} - 2 \, q^{5} + q^{4} + 3 \, q^{2} = q^{10} - 2 \, q^{7} - q^{5} + q^{4} + q^{2}
\end{align*}
implying that $-q^5 + 2q^2 = 0$. Therefore, such $q$ also cannot exist.
\begin{thm}
	The Alternating forms graph $\alt{5,q}$ is not triply transitive for any prime power $q$.\label{thm:alternating-form}
\end{thm}
\subsection{Bilinear forms graphs}
The bilinear forms graph $H_q(2,e)$ is the graph whose vertex set consists of the $2\times e$ matrices over $\mathbb{F}_q$, where two matrices are adjacent if their difference has rank equal to $1$. The bilinear forms graph $H_q(2,e)$ is a rank $3$ graph with parameters
\begin{align*}
	\left(q^{2e},(q+1)(q^e-1),q^e+(q-2)(q+1),q(q+1)\right).
\end{align*}
We claim that $H_q(2,e)$ is not triply regular. The graph $H_q(2,e)$ is clearly a Cayley graph on the additive group of $2\times e$ matrices over $\mathbb{F}_q$, and with connection set consisting of all rank $1$ matrices. Now, consider the matrices
\begin{align*}
	A &= 
	\begin{bmatrix}
		1 & 0 &0& \ldots & 0\\
		0 & 0 &0& \ldots & 0
	\end{bmatrix}\\
	B &= 
	\begin{bmatrix}
		0 & 1 &0& \ldots & 0\\
		0 & 0 &0& \ldots & 0
	\end{bmatrix}\\
	C&=
	\begin{bmatrix}
		0 & 0 &0& \ldots & 0\\
		1 & 0 &0& \ldots & 0
	\end{bmatrix}.
\end{align*}
It is clear that $A,B$, and $C$ are vertices of the first subconstituent of $H_q(2,e)$ with respect to the zero matrix. Moreover, $A$ is adjacent to both $B$ and $C$. In the first subconstituent, we will show that the number of common neighbours of $A$ and $B$ and the one of $A$ and $C$ are not equal.
\begin{prop}
	The number of common neighbours of $A$ and $B$ is $q^e-3$.\label{prop:bil-form-1}
\end{prop}
\begin{proof}
	Let 
	\begin{align*}
		X &=
		\begin{bmatrix}
			a_1 & a_2 & a_3& \ldots & a_e\\
			b_1 & b_2 & b_3& \ldots & b_e
		\end{bmatrix}
	\end{align*}
	be a common neighbour of $A$ and $B$ in the first subconstituent. Since $X$ is adjacent to $A$, we know that
	\begin{align}
		X-A = 
		\begin{bmatrix}
			a_1-1 & a_2 & a_3& \ldots & a_e\\
			b_1 & b_2 & b_3& \ldots & b_e
		\end{bmatrix}\label{eq:X-A}
	\end{align}
	and 
	\begin{align}
		X-B = 
		\begin{bmatrix}
			a_1 & a_2-1 & a_3& \ldots & a_e\\
			b_1 & b_2 & b_3& \ldots & b_e
		\end{bmatrix}\label{eq:X-B}
	\end{align}
	are both rank $1$.
	
	\begin{itemize}
		\item If $a_1 = 1$, then we consider the two cases whether $b_1$ is zero or not. If $b_1 = 0$, then \eqref{eq:X-B} becomes	
		\begin{align*}
			X-B = 
			\begin{bmatrix}
				1 & a_2-1 & a_3& \ldots & a_e\\
				0 & b_2 & b_3& \ldots & b_e
			\end{bmatrix}.
		\end{align*}
		As $X-B$ has rank $1$, we deduce that $b_2 = b_3 = \ldots = b_e = 0$. Consequently, we have
		\begin{align*}
			X &=
			\begin{bmatrix}
				1 & a_2 & a_3& \ldots & a_e\\
				0 & 0 & 0& \ldots & 0
			\end{bmatrix}.
		\end{align*}
		It is immediate to see that there are exactly 
		\begin{align}
			q^{e-1}-1\label{eq:common-neighbours-1}
		\end{align}
		 such matrices since $X\neq A$.
		
		Next, we consider the case where $b_1\neq 0$. In this case, \eqref{eq:X-A} becomes
		\begin{align*}
			X-A = 
			\begin{bmatrix}
				0 & a_2 & a_3& \ldots & a_e\\
				b_1 & b_2 & b_3& \ldots & b_e
			\end{bmatrix}.
		\end{align*}
		It is therefore immediate that $a_2=a_3=\ldots = a_e = 0$. Therefore, we have
		\begin{align*}
			X &=
			\begin{bmatrix}
				1 & 0 & 0& \ldots & 0\\
				b_1 & b_2 & b_3& \ldots & b_e
			\end{bmatrix}.
		\end{align*}
		Using the fact that $X$ also has rank $1$, we must have $b_2 = b_3 =\ldots = b_e = 0$, and so
		\begin{align*}
			X &=
			\begin{bmatrix}
				1 & 0 & 0& \ldots & 0\\
				b_1 & 0 & 0& \ldots & 0
			\end{bmatrix}.
		\end{align*}
		As $X-B$ has rank $1$, we must have $b_1 = 0$, which is a contradiction.
		\item Now, suppose that $a_1\neq 1$. In this case, $a_1-1 \neq 0$, and so there exists $k\in \mathbb{F}_q$ such that $b_1 = k(a_1-1)$ and $b_i = k a_i$ for $2\leq i\leq e$. Hence, \eqref{eq:X-A} becomes
		\begin{align*}
			X-A = 
			\begin{bmatrix}
				a_1-1 & a_2 & a_3& \ldots & a_e\\
				k(a_1-1) & ka_2 & ka_3& \ldots & ka_e
			\end{bmatrix},
		\end{align*} 
		and \eqref{eq:X-B} becomes
		\begin{align}
			X-B = 
			\begin{bmatrix}
				a_1 & a_2-1 & a_3& \ldots & a_e\\
				k(a_1-1) & ka_2 & ka_3& \ldots & ka_e
			\end{bmatrix}.\label{eq:X-B-2}
		\end{align}
		If $k = 0$, then $X-A$ and $X-B$ both have rank $1$, and we have
		\begin{align*}
			X &=
			\begin{bmatrix}
				a_1 & a_2 & a_3& \ldots & a_e\\
				0 & 0 & 0& \ldots & 0
			\end{bmatrix}.
		\end{align*}
		By noting that $a_1\neq 1$ and $X\notin \{0,A,B\}$, there are exactly 
		\begin{align}
			(q^{e}-1)-q^{e-1}-1 = q^e-q^{e-1}-2\label{eq:common-neighbours-2}
		\end{align}
		 such elements if $k = 0$.
		
		If $k\neq 0$, then $k(a_1-1)\neq 0$, and using the fact that the matrix $X-B$ in \eqref{eq:X-B-2} has rank $1$, we deduce that there exists $k^\prime \in \mathbb{F}_q$ such that $a_1= k^\prime k(a_1-1), a_2-1 = k^\prime ka_2$, and $a_i = k^\prime k a_i$, for all $3\leq i\leq e$. From the latter, we have $k^\prime = k^{-1}$ or  $a_i = 0$ for all $3\leq i\leq e$.
		
		\noindent{\bf Case~1:} $k^\prime = k^{-1}$.
		
		In this case, we have $a_1 = a_1-1$, and $ a_2-1 = a_2$, which are impossible.
		
		\noindent{\bf Case~2:} Assume that $k^\prime k\neq 1$ and $a_i = 0$ for all $3\leq i\leq e$.
		
		In this case, since $a_1= k^\prime k(a_1-1), a_2-1 = k^\prime ka_2$, we know that
		\begin{align*}
			a_1 = \frac{1}{k^\prime k-1} \mbox{ and } a_2 = -\frac{1}{k^\prime k-1}.
		\end{align*}
		We conclude that 
		\begin{align*}
			X &=
			\begin{bmatrix}
				\frac{1}{k^\prime k-1} & -\frac{1}{k^\prime k-1} & 0& \ldots & 0\\
				\frac{k(2-k^\prime k)}{k^\prime k-1} & -\frac{k}{k^\prime k-1} & 0& \ldots & 0
			\end{bmatrix}.
		\end{align*}
		For the above matrix to be of rank $1$, we must have $-k+k(2-k^\prime k) = 0$. Since $k\neq 0$, we must have $1 = kk^\prime$, which is impossible.
	We conclude that $X$ is not a common neighbour of $A$ and $B$ when $a_1\neq 0$.
	\end{itemize}
	Hence, by combining \eqref{eq:common-neighbours-1} and \eqref{eq:common-neighbours-2}, $A$ and $B$ have $q^{e-1}-1+q^e-q^{e-1}-2 = q^e-3$ common neighbours.
\end{proof}

Next, we determine the number of common neighbours of $A$ and $C$.
\begin{prop}
	The number of common neighbours of $A$ and $C$ is $q^2-3$.\label{prop:bil-form-2}
\end{prop}
\begin{proof}
	Again, let 
	\begin{align*}
		X &=
		\begin{bmatrix}
			a_1 & a_2 & a_3& \ldots & a_e\\
			b_1 & b_2 & b_3& \ldots & b_e
		\end{bmatrix}
	\end{align*}
	be a common neighbour of $A$ and $C$ in the first subconstituent. We will distinguish the cases whether $a_1$ is zero or not.
	\begin{itemize}
		
		\item Suppose that $a_1 = 0$. Then, we have
		\begin{align*}
			X &=
			\begin{bmatrix}
				0 & a_2 & a_3& \ldots & a_e\\
				b_1 & b_2 & b_3& \ldots & b_e
			\end{bmatrix}.
		\end{align*}
		
		If $b_1\neq 0$, then $a_2 = a_3 = \ldots = a_e = 0$ since $X$ has rank $1$. Since $X-A$ also has rank $1$, we must have $b_2 = b_3 = \ldots = b_e = 0$. Hence, 
		\begin{align*}
			X &=
			\begin{bmatrix}
				0 & 0 & 0& \ldots & 0\\
				b_1 & 0 & 0& \ldots & 0
			\end{bmatrix}.
		\end{align*}
		As, $X-C$ has rank $1$, we must have $b_1 \neq 1$. There are 
		\begin{align}
			q-2\label{eq:number3}
		\end{align}
		such matrices.
		
		If $b_1 = 0$, then we have
		\begin{align*}
			X &=
			\begin{bmatrix}
				0 & a_2 & a_3& \ldots & a_e\\
				0 & b_2 & b_3& \ldots & b_e
			\end{bmatrix}.
		\end{align*}
		As $X-A$ has rank $1$, we have $b_2 = b_3 = \ldots = b_e = 0$. Similarly, since $X-C$ has rank $1$, we have $a_2 = a_3= \ldots = a_e=0$. In this case, $X$ is the zero matrix, which does not belong to the first subconstituent of $H_q(2,e)$. So there are no common neighbour $X$ where $b_1 = 0$.
		
		\item Suppose that $a_1\neq 0$. Then, there exists $k\in \mathbb{F}_q$ such that $b_i = ka_i$ for all $1\leq i\leq e$. Then, $X$ becomes
		\begin{align*}
			X &=
			\begin{bmatrix}
				a_1 & a_2 & a_3& \ldots & a_e\\
				ka_1 & ka_2 & ka_3& \ldots & ka_e
			\end{bmatrix}.
		\end{align*}
		Moreover, we have
		\begin{align}
			X -A=
			\begin{bmatrix}
				a_1-1 & a_2 & a_3& \ldots & a_e\\
				ka_1 & ka_2 & ka_3& \ldots & ka_e
			\end{bmatrix}\label{eq:X-A-case2}
		\end{align}
		and
		\begin{align}
			X -C=
			\begin{bmatrix}
				a_1 & a_2 & a_3& \ldots & a_e\\
				ka_1-1 & ka_2 & ka_3& \ldots & ka_e
			\end{bmatrix}\label{eq:X-C}
		\end{align} 
		
		If $k = 0$, then it is clear that $X-A$ has rank $1$ unless $a_1=1$ and $a_2 = a_3 = \ldots = a_e = 0$. In addition, $X-C$ has rank $1$ only when $a_2 = a_3 = \ldots = a_e = 0.$ Therefore, 
		\begin{align*}
			X &=
			\begin{bmatrix}
				a_1 & 0 & 0& \ldots & 0\\
				0 & 0 & 0& \ldots & 0
			\end{bmatrix},
		\end{align*}
		where $a_1\not\in \{0,1\}$. Hence, there are exactly
		\begin{align}
			q-2\label{eq:number1}
		\end{align}
		such matrices $X$ whenever $k = 0$.
		
		Assume now that $k\neq 0$. Then $ka_1\neq 0$. The first row of the matrix in \eqref{eq:X-A-case2} is a scalar multiple of the second row. Therefore, there exists $k^\prime \in \mathbb{F}_q$ such that $a_1-1 = k^\prime k a_1$ and $a_i = k^\prime k a_i$ for all $2 \leq i\leq e$. 
		
		Since $a_1-1 = k^\prime k a_1$, it is clear that $k^\prime k \neq 1$ and $a_1 = \tfrac{1}{1-k^\prime k}$. Using the fact that $a_i = k^\prime k a_i$ for all $2 \leq i\leq e$, we must have $a_i = 0$ for $2\leq i\leq e$. Therefore, we have
		\begin{align*}
			X &=
			\begin{bmatrix}
				a_1 & 0 & 0& \ldots & 0\\
				ka_1 & 0 & 0& \ldots & 0
			\end{bmatrix}
			=
			\begin{bmatrix}
				\tfrac{1}{1-k^\prime k} & 0 & 0& \ldots & 0\\
				\tfrac{k}{1-k^\prime k} & 0 & 0& \ldots & 0
			\end{bmatrix}.
		\end{align*}
		Since $k\neq 0$ and $k^\prime \neq k^{-1}$, there are exactly
		\begin{align}
			(q-1)(q-1) = q^2-2q+1\label{eq:number2}
		\end{align}
		matrices of this form.

	\end{itemize}
	
	We conclude from \eqref{eq:number3}, \eqref{eq:number1}, and \eqref{eq:number2} that there are exactly $(q^2-2q+1)+(q-2)+(q-2) = q^2-3$ common neighbours of $A$ and $C$ in the first subconstituent.
\end{proof}

From Proposition~\ref{prop:bil-form-1} and Proposition~\ref{prop:bil-form-2}, the bilinear forms graph $H_q(2,e)$ is not triply regular for $e\geq 3$. Hence, $H_q(2,e)$ is not triply transitive for $e\geq 3$. The graph $H_q(2,2)$ is isomorphic to the affine polar graph $VO^+_4(q)$, so by \S\ref{subs:affine-polar-graph}, we conclude that $H_q(2,2)$ is not triply transitive unless $q = 2$. We state all of these in the next theorem.
\begin{thm}
	The bilinear forms graph $H_q(2,e)$ is triply transitive if and only if $q = e = 2$.\label{thm:bilinear-form}
\end{thm}

\subsection{The graph $NU_m(2)$}
Let $V$ be a vector space of dimension $m\geq 3$ over $\mathbb{F}_4$ equipped with a non-degenerate Hermitian form $f$. The graph $NU_m(2)$ is the graph on non-isotropic points of $V$, adjacent when they are joined by a tangent. The graph $NU_m(2)$ is a rank $3$ graph, and its complement $\overline{NU_m(2)}$ is the graph on non-isotropic points, where two vertices are adjacent when they are orthogonal. The graph $\overline{NU_3(2)}$ is disconnected, and thus imprimitive. The graph $\overline{NU_m(2)}$, for $m\geq 4$, is connected and is locally $\overline{NU_{m-1}(2)}$, so its first subconstituent is strongly regular.

The graph $\overline{NU_m(2)}$ is a rank $3$ graph with $\tfrac{2^{m-1}(2^m-\varepsilon)}{3}$ vertices, where $\varepsilon = (-1)^m$. If $m$ is even, then $\tfrac{2^{m-1}(2^m-\varepsilon)}{3} = \tfrac{2^{m-1}(2^m-1)}{3}$ cannot be a square since $2^m-1$ is odd and ${m-1}$ is odd. Similarly, if $m$ is odd, then $\tfrac{2^{m-1}(2^m-\varepsilon)}{3} = \tfrac{2^{m-1}(2^m+1)}{3}$ is not a square. Therefore, $\overline{NU_m(2)}$ does not have Latin square or negative Latin square parameters. As the Krein parameters $\mathrm{q}_{11}^1$ and $\mathrm{q}_{22}^2$ are both positive, by Lemma~\ref{lem:srg-not-triply-regular}, we deduce that the graph $\overline{NU_m(2)}$ is not triply transitive. 

\begin{thm}
	The graph $NU_m(2)$ is not triply transitive for any $m\geq 3$.\label{thm:non-singular-unitary}
\end{thm}

\subsection{The graph $\no^\varepsilon_{2m}(q)$}
For $m\geq 3$ and $\varepsilon = \pm 1$, the graphs $\no^\varepsilon_{2m}(q)$ for $q \in \{2,3\}$ are rank $3$ graphs. See \cite[\S3.1.2]{brouwer2022strongly} for details on the construction. 

The number of vertices of $\no^\varepsilon_{2m}(2)$ is equal to $2^{2m-1}-\varepsilon 2^{m-1}$. If $\varepsilon = 1$, then there are $2^{2m-1}-2^{m-1} = 2^{m-1}(2^m-1)$. Hence, if $m$ is even, then $2^{m-1}(2^m-1)$ is not a square. If $m$ is odd, then $m-1$ is even and $2^m-1$ is a square if and only if $m = 1$. Consequently, $2^{m-1}(2^m-1)$ is not a square for $m\geq 2$. Hence, $\no^\varepsilon_{2m}(2)$ does not have a Latin square or a negative Latin square parameters. Moreover, the Krein parameters $\mathrm{q}_{11}^1$ and $\mathrm{q}_{22}^2$ of $\no^\varepsilon_{2m}(2)$ are both positive, so we may apply Lemma~\ref{lem:srg-not-triply-regular}.
Therefore, the graph $\no^\varepsilon_{2m}(2)$ is not triply transitive.

Similarly, the number of vertices of $\no^\varepsilon_{2m}(3)$ is $\tfrac{3^{m-1}(3^m-\varepsilon)}{2}$. If $\varepsilon = +1$ and $m$ is even, then $\frac{3^{m-1}(3^m-1)}{2}$ is not a square. If $\varepsilon = +1$ and $m$ is odd, then $m-1$ is even and it is not hard to show that
\begin{align*}
	\frac{(3^m-1)}{2}
\end{align*}
is not a square unless $m\in \{1,2\}$. Therefore, we conclude that  $\no^+_{2m}(3)$ does not have Latin square or negative Latin square parameters.

If $\varepsilon = -1$, then $\tfrac{3^{m-1}(3^m-\varepsilon)}{2} = \tfrac{3^{m-1}(3^m+1)}{2}$ which is clearly not a square if $m$ is even. If $m$ is odd, then $m-1$ is even and it is not hard to see that
\begin{align*}
	\frac{(3^m+1)}{2}
\end{align*}
is a square only when $m = 0$. We conclude that the number of vertices of $\no^-_{2m}(3)$ is never a square. Consequently, $\no^-_{2m}(3)$ does not have Latin square or negative Latin square parameters.

The Krein parameters $\mathrm{q}_{11}^1$ and $\mathrm{q}_{22}^2$ are both positive for the graph $\no_{2m}^\varepsilon(3)$, so by Lemma~\ref{lem:srg-not-triply-regular}, we deduce that $\no_{2m}^\varepsilon(3)$ is also not triply transitive for $m\geq 1$ and $\varepsilon = \pm 1.$

\begin{thm}
	The graph $\no_{2m}^\varepsilon(q)$ is not triply transitive for any $q\in \{2,3\}$, $\varepsilon = \pm 1$, and $m\geq 3$.\label{thm:non-singular-orthogonal-even}
\end{thm}

\subsection{The graphs $\no^\varepsilon_{2m+1}(q)$}

Let $m\geq 1$ and $\varepsilon = \pm 1$ such that $(m,\varepsilon) \neq (1,-1)$, $q \in \{3,4,8\}$, and consider the rank $3$ graph $\no^\varepsilon_{2m+1}(q)$ defined in \cite[\S3.1.4]{brouwer2022strongly}. The number of vertices of $\no^\varepsilon_{2m+1}(q)$ is $\binom{q^n}{2}$ when $\varepsilon = -1$ and $\binom{q^n+1}{2}$, when $\varepsilon = 1$. As these numbers are non-squares for $(q,m,\varepsilon) \neq (3,2,-1)$ and $(q,m,\varepsilon)=(8,1,+1)$, the graph $\no^\varepsilon_{2m+1}(q)$ does not have Latin square or a negative Latin square parameters. The graph $\no_{3}^+(8)$ has parameters $(36, 14, 7, 4)$, so it does not have Latin square or negative Latin square parameters. The graph $\no_{5}^-(3)$ has parameters $(36, 20, 10, 12)$, which are Latin square parameters. For the latter, we use \verb*|Sagemath| to show that $\dim(\tilde{T}) = 16$, so it cannot be triply transitive. For the remaining cases that do not have Latin square or negative Latin square parameters, the Krein parameters $\mathrm{q}_{11}^1$ and $\mathrm{q}_{22}^2$ of $\no^\varepsilon_{2m+1}(q)$ are both positive, so by Lemma~\ref{lem:srg-not-triply-regular}, the graph $\no^\varepsilon_{2m+1}(q)$ is not triply transitive.

\begin{thm}
	The graph $\no_{2m+1}^\varepsilon(q)$ is not triply transitive for any $q\in \{3,4,8\}$, $\varepsilon = \pm 1$, and $m\geq 1$.\label{thm:non-singular-orthogonal-odd}
\end{thm}

\subsection{Distance-$2$ graphs of polar spaces of rank $5$}
Let $\mathcal{P}$ be a classical polar space of rank $5$ and order $(q,1)$. Consider the distance-$2$ graph $\Gamma$ of the dual polar graph arising from $\mathcal{P}$. This graph is a rank $3$ graph with parameters $(n,k,\lambda,\mu)$ given in \cite[Theorem~2.2.20]{brouwer1990graphs}. In particular, the negative eigenvalue of $\Gamma$ is $-q^2-1$. If $\Gamma$ has Latin square or negative Latin square parameters, then by Lemma~\ref{lem:rank3-LS-nLS}, we deduce that the parameters of $\Gamma$ are of the form $(4t^2,t(2t\pm1),t(t\pm1),t(t\pm1))$, where $t = q^2+1$. This is however impossible since the number of vertices of $\Gamma$ is $(q^4+1)(q^3+1)(q^2+1)(q+1)$.

Therefore, the graph $\Gamma$ does not have Latin square or negative Latin square parameters. The Krein parameters $\mathrm{q}_{11}^1$ and $\mathrm{q}_{22}^2$ are both positive, so by Lemma~\ref{lem:srg-not-triply-regular}, we deduce that the rank $3$ graph arising from this construction cannot be triply transitive.

\begin{thm}
	If $\Gamma$ is the distance-2 graph of a classical polar space of rank $5$, then $\Gamma$ is not triply transitive.\label{thm:rank5}
\end{thm}

\subsection{The van Lint-Schrijver graphs}

Let $q = p^f$, where $p$ is a prime number, and let $e\mid (q-1)$ be an odd prime such that $p$ is primitive $\mod{e}$ such that $f = (e-1)\ell$, for some integer $t\geq 1$. Assume that $q = em+1$. Fix a primitive root of unity $\mathbb{F}_q$. Let $K\leq \mathbb{F}_q^*$ be the subgroup of the $e$-th powers of elements of $\mathbb{F}_q^*$. The van Lint-Schrijver graph $\Gamma_q$ is defined as the Cayley graph of the additive group $\mathbb{F}_q$ with connection set equal to $K$. Hence, $\Gamma_q$ has parameters
\begin{align*}
	\left(q,\tfrac{q-1}{e},\lambda,\mu\right)
\end{align*} 
where $\lambda$ and $\mu$ are given in \cite{eindhoven1981construction}. Its complement therefore has parameters
\begin{align*}
	\left(q,\tfrac{(e-1)(q-1)}{e},q-\tfrac{2(q-1)}{e}+\mu-2,q-\tfrac{2(q-1)}{e}+\lambda\right).
\end{align*}

By \cite{eindhoven1981construction}, the graph $\Gamma_q$ has Latin square parameters if $\ell$ is odd, and negative Latin square parameters if $\ell$ is even. Using Lemma~\ref{lem:rank3-LS-nLS}, if $\Gamma_q$ is not triply transitive unless $\Gamma_q$ or $\overline{\Gamma}_q$ has Latin square parameters $(4t^2,t(2t-1),t(t-1),t(t-1))$ when $\ell$ is odd and negative Latin square parameters $(4t^2,t(2t+1),t(t+1),t(t+1))$ when $\ell$ is even. 

Assume first that $\Gamma_q$ has these parameters. Since the valency of $\Gamma_q$ is $\tfrac{q-1}{e}$ and $4t^2 = q$, we must have $e = \tfrac{2t\pm 1}{t}$ which can only be an integer when $t = 1$. In both cases, $\Gamma_q$ has parameters $(4,3,2,2)$, which is impossible, or $(4,1,0,0)$, which is imprimitive. 

Similarly, if $\overline{\Gamma}_q$ has these Latin square or negative Latin square parameters, then it is not hard to show that
\begin{align*}
	e = \frac{2t\pm 1}{t\pm 1}
\end{align*}
which cannot be a prime number.

Hence, we deduce the following result.
\begin{thm}
	The van Lint-Schrijver graphs are not triply transitive.\label{thm:vanLint-Schrijver}
\end{thm}

\section{Proof of the main results}\label{Sec:main-results}

\subsection{Triple transitivity in rank $3$ graphs}\label{subsect:non-triply-transitive}
Let us first identify the graphs that are not triply transitive.
\begin{enumerate}[(1)]
	\item Except for the Higman-Sims graph and McLaughlin graph, no other graphs in Table~\ref{tab} are triply transitive. The Higman-Sims graph has negative Latin square parameters. The McLaughlin graph has neither Latin square nor negative Latin square parameters.
	\item The Johnson graphs $J(n,2)$, where $n\geq 5$ (Theorem~\ref{thm:triangular}). 
	\item The Grassmann graph $J_q(n,2)$, where $n\geq 4$ (Theorem~\ref{thm:Grassmann}).
	\item The Paley graph $P(q)$, for $q>9$, (Theorem~\ref{thm:Paley}).
	\item The Peisert graph $P^*(p^{2t})$, for $(p,t)\neq (3,2)$, (Theorem~\ref{thm:Peisert}).
	\item The $E_{6,1}(q)$ graphs (Theorem~\ref{thm:e61}).
	\item The collinearity graphs of classical polar spaces $\mathcal{P}$ of rank $n\geq 2$, except when $\mathcal{P} \in \{ O_6^-(q),O_4^+(q) \}$ (see Theorems~\ref{thm:polar-graph-rank-3},\ref{thm:polar-graph-rank-2}).
	\item The collinearity graphs of the dual of classical polar spaces $\mathcal{P}$ of rank $2$ except when $\mathcal{P} = O^+_4(q)$ (see Theorem~\ref{thm:polar-dual-graph-rank-2}).
	\item The affine polar graphs $\vo_{2m}^\varepsilon(q)$, where $q>2$ and $m\geq 2$ (see Theorem~\ref{thm:affine-polar-graph}).
	\item The affine half spin $\vd_{5,5}(q)$, where any prime power $q$ (see Theorem~\ref{thm:affine-half-spin}).
	\item The graph $\vsz(q)$ for any prime power $q = 2^{2e-1}$, where $e\geq 2$ (see Theorem~\ref{thm:ST}).
	\item The alternating forms graphs $\alt{5,q}$ for any prime power $q$ (see Theorem~\ref{thm:alternating-form})
	\item The bilinear forms graphs $H_q(2,e)$ for any prime power $q$ and integer $e\geq 2$ such that $(q,e)\neq (2,2)$ (see Theorem~\ref{thm:bilinear-form}).
	\item The graph $NU_{m}(2)$ for any $m\geq 3$ (see Theorem~\ref{thm:non-singular-unitary}).
	\item The graph $\no_{2m}^\varepsilon(q)$ for $m\geq 3$, $\varepsilon = \pm 1$, and $q\in \{2,3\}$ (see Theorem~\ref{thm:non-singular-orthogonal-even}).
	\item The graph $\no_{2m+1}^\varepsilon(q)$ for $m\geq 3$, $\varepsilon = \pm 1$, and $q\in \{3,4,8\}$ (see Theorem~\ref{thm:non-singular-orthogonal-odd}).
	\item The distance-2 graph of polar spaces of rank $5$ (see Theorem~\ref{thm:rank5}).
	\item The van Lint-Schrijver graph (see Theorem~\ref{thm:vanLint-Schrijver}).
\end{enumerate}

Next, we identify the graphs that are triply transitive. If $\Gamma$ is triply transitive, then $\Gamma$ is among the following graphs.
\begin{enumerate}[(a)]
	\item The Higman-Sims graph.\label{HS2}
	\item  the McLaughlin graph.\label{McL}
	\item A complete multipartite graph with $n$ parts of size $m$ in Theorem~\ref{thm:imprimitive}.\label{imprimitive2}
	\item The Paley graphs $P(q)$, where $q\in \{5,9\}$ in Theorem~\ref{thm:Paley}. Note that $P(9)$ and $P^*(9)$ are isomorphic.\label{Paley2}
	\item The collinearity graph of $O_6^-(q)$.\label{O6}
	\item The affine polar graph $\vo_{2m}^\varepsilon(2)$.\label{affine-polar-graph}
	\item The $n\times n$ grid.\label{grid2}
\end{enumerate}

\subsection{Proof of Theorem~\ref{thm:main}}
Let $\Gamma$ be a strongly-regular graph as in Hypothesis~\ref{hyp0}. Assume that $\Gamma$ does not have a Latin square or negative Latin square parameters. If $\Gamma$ is triply transitive, then it is one of the graphs in \eqref{HS2}-\eqref{grid2}. However, since $\Gamma$ does not have a Latin square or negative Latin square parameters, $\Gamma$ cannot be the graphs in \eqref{HS2}, \eqref{imprimitive2} with $n = m$, \eqref{Paley2} where $q = 9$, \eqref{affine-polar-graph} by Remark~\ref{rmk:affine-polar-graph}, or \eqref{grid2}.

Hence, $\Gamma$ is the McLaughlin graph, a complete multipartite graph with $n$ parts of size $m$ where $n\neq m$, a $5$-cycle, or the collinearity graph of $O_6^-(q)$. These are the graphs listed in the statement of the theorem.
Except possibly for the collinearity graph of $O_6^-(q)$, all graphs in the list are triply transitive by Row~18 in Table~\ref{tab}, Theorem~\ref{thm:imprimitive} and Theorem~\ref{thm:Paley}. This completes the proof.

\subsection{Proof of Theorem~\ref{thm:main2}}

Let $\Gamma$ be a graph as in Hypothesis~\ref{hyp0} which has Latin square or negative Latin square parameters. If $\Gamma$ is triply transitive, then $\Gamma$ is one of the graphs in \eqref{HS2}-\eqref{grid2} in Section~\ref{subsect:non-triply-transitive}. Hence, $\Gamma$ is one of the graphs in \eqref{HS2}, \eqref{imprimitive2} with $n = m$, \eqref{Paley2} with $q = 9$, \eqref{affine-polar-graph}, or \eqref{grid2}.
These are exactly the graphs in the statement of Theorem~\ref{thm:main2}. Except for the graph in \eqref{affine-polar-graph}, all graphs in this list are triply transitive\footnote{see Theorems~\ref{thm:imprimitive},~\ref{thm:Peisert},~\ref{thm:grid},~\ref{thm:polar-graph-rank-2},~\ref{thm:affine-polar-graph}}. This completes the proof of Theorem~\ref{thm:main2}. 
\section{Concluding remarks and open questions}\label{Sec:conclusion}

In this paper, we studied the strongly-regular graphs that are triply transitive. In Theorem~\ref{thm:main} and Theorem~\ref{thm:main2}, we have identified all possible triply transitive graphs, and showed that the majority of the graphs in this list are triply transitive. In order to obtain a full classification, one would need to determine whether or not the graphs in Theorem~\ref{thm:main}~\eqref{o} and Theorem~\ref{thm:main2}~\eqref{vo} are triply transitive. We have checked using \verb*|Sagemath| that these graphs are triply transitive for small values of their parameters. Therefore, we are inclined to believe that these graphs are triply transitive, and we make the following conjectures.


\begin{qst}
	 Is the affine polar graph $VO^\varepsilon_{2m}(2)$ triply transitive, for $m\geq 1$ and $\varepsilon = \pm 1$?\label{conj1}
\end{qst}

\begin{qst}
	 Is the collinearity graph of polar space $O_6^-(q)$ triply transitive, when $q$ is a prime power?\label{conj2}
\end{qst}

If the answer to the above questions is affirmative, then the following conjecture (a complete classification) follows as a corollary.
\begin{conj}
	Let $\Gamma$ be a graph as in Hypothesis~\ref{hyp0}. The graph $\Gamma$ is triply transitive if and only if $\Gamma$ is one of the following graphs.
	\begin{enumerate}[(a)]
		\item A complete multipartite graph with $n$ parts of size $m$.
		\item The cycle on $5$ vertices.
		\item The McLaughlin graph.
		\item The Higman-Sims graph.
		\item The Peisert graph $P^*(9)$ (which is isomorphic to the Paley graph on $9$ vertices).
		\item The $n\times n$ grid for some integer $n\geq 2$.
		\item The collinearity graph of the polar space  $O_6^-(q)$.
		\item The affine polar graph $VO^\varepsilon_{2m}(2)$, for $m\geq 2$ and $\varepsilon =\pm 1$.
	\end{enumerate}
\end{conj}

In Theorem~\ref{thm:decomposition-T-tilde}, we showed that the block dimension decomposition of the centralizer algebra of the vertex stabilizer of the full automorphism group of a rank $3$ graph $\Gamma$ is of the form
\begin{align}
	\begin{bmatrix}
		1 & 1 & 1\\
		1 & r_1 & t\\
		1 & t & r_2
	\end{bmatrix}\label{eq:decomposition}
\end{align}
where $t\geq 1$ is an integer, $r_1$ and $r_2$ are the rank of the vertex stabilizer in its actions on the first and second subconstituents, respectively. If $\Gamma$ is a Paley graph on a prime number of vertices $p$, then the group $\Aut{\Gamma}$ is a Frobenius group, and the block dimension decomposition of the $\tilde{T}$ is of the form
\begin{align*}
	\begin{bmatrix}
		1 & 1 & 1\\
		1 & \tfrac{p-1}{2} & \tfrac{p-1}{2}\\
		1 & \tfrac{p-1}{2} & \tfrac{p-1}{2}
	\end{bmatrix}.
\end{align*}
Hence, $r_1 = r_2 = t = \tfrac{p-1}{2}$. If $\Gamma $ is the Peisert graph $ P^*(p^2)$, then the full automorphism group is also a Frobenius group, and $r_1 = r_2 = t = \tfrac{p^2-1}{2}$.
\begin{qst}
	Except for the Paley and Peisert graphs, are there other rank $3$ graphs where $r_1 = r_2 = t$ in \eqref{eq:decomposition}?
\end{qst}

We have observed in our computations with \verb*|Sagemath| that $t< \max(r_1,r_2)$ (see \eqref{eq:decomposition}) in all the examples. Hence, we ask the following equivalent question.
\begin{qst}
	Does the inequality $t< \max(r_1,r_2)$ always hold, where $r_1,r_2$ and $t$ are given in \eqref{eq:decomposition} for rank $3$ graphs?
\end{qst}

If the answer to the above question is affirmative, then Question~\ref{conj1} and Question~\ref{conj2} both hold, and hence we obtain a full classification of triply transitive graphs.

%

\end{document}